\documentclass[11pt,a4paper]{article}
\usepackage[margin=0.75in]{geometry}
\usepackage[utf8]{inputenc}
\usepackage{cmap}
\usepackage[T1]{fontenc}
\usepackage{lmodern}
\usepackage[square, numbers]{natbib}     
\usepackage[colorlinks=true]{hyperref}
\usepackage{amssymb, amsthm, amsmath, amsfonts}
\usepackage{authblk, dsfont, mathrsfs, color, bm, enumitem, mathtools}
\usepackage{subfigure, graphicx, transparent}
\usepackage{dirtytalk}
\graphicspath{{figures/}}
\usepackage[dvipsnames]{xcolor}
\usepackage{comment}
\hypersetup{colorlinks, linkcolor={red!80!black}, citecolor={blue!80!black}, urlcolor={ForestGreen}}
\usepackage{algorithm}
\usepackage{algpseudocode}

\theoremstyle{definition}
\newtheorem{theorem}{Theorem}[section]
\newtheorem{definition}[theorem]{Definition}

\newtheorem{proposition}[theorem]{Proposition}

\newtheorem{lemma}[theorem]{Lemma}
\newtheorem{remark}[theorem]{Remark}
\newtheorem{example}[theorem]{Example}

\DeclareUnicodeCharacter{211D}{\ensuremath{\mathbb R}}

\DeclarePairedDelimiterX{\inner}[2]{\langle}{\rangle}{#1, #2}

\title{Diffusion-Robust Optimization over Graphs}
\author[1]{Liviu Aolaritei$^*$}
\author[2]{Ricky Huang$^*$}
\author[1,3,4]{Michael I. Jordan}
\author[2]{Paul Grigas}

\affil[1]{Department of Electrical Engineering and Computer Sciences, UC Berkeley, USA}
\affil[2]{Department of Industrial Engineering and Operations Research, UC Berkeley, USA}
\affil[3]{Department of Statistics, UC Berkeley, USA}
\affil[4]{Inria Paris, France \protect\\
\texttt{\{liviu.aolaritei,\,yxhuang,\,pgrigas\}@berkeley.edu, jordan@cs.berkeley.edu}}

\date{}

\begin{document}
\maketitle

\begin{abstract}
We introduce a diffusion-based uncertainty model for robust optimization on directed graphs, in which perturbations of edge weights propagate along adjacent edges and satisfy conservation constraints at nodes. This topology-aware structure is natural in networked systems where uncertainty is induced by flows and local interactions, including transportation, logistics, communication, and energy networks. We analyze how such diffusive uncertainty reshapes the computational landscape of robust graph optimization. We focus on two canonical combinatorial graph problems, shortest path and the traveling salesman problem (TSP), which provide complementary benchmarks: shortest path is polynomial-time solvable in the nominal setting, whereas TSP is already \textsf{NP}-hard. We show that, for shortest path, propagation depth induces a sharp transition between tractable and intractable robust counterparts. For the traveling salesman problem, robustness often adds no computational complexity beyond ordinary TSP, because the structure of Hamiltonian cycles makes the fixed-tour adversarial problem collapse to explicit formulas. Together, these results show that topology-aware uncertainty can fundamentally change robust combinatorial optimization, with tractability governed by the interaction between propagation, budget geometry, and the structure of feasible solutions.
\end{abstract}

\begingroup
\renewcommand\thefootnote{}
\footnotetext{$^*$: Equal contribution.}
\endgroup

\section{Introduction}
\label{sec:intro}

Optimization on networks lies at the heart of discrete optimization. A directed graph is a deceptively simple mathematical structure, yet it supports some of the deepest results in algorithm design and computational complexity. On this structure, shortest paths, minimum cuts, and network flow problems admit elegant polynomial-time algorithms built on refined combinatorial insights, while the traveling salesman problem and many network design problems remain \textsf{NP}-hard and define canonical boundaries of tractability \cite{ahuja1988network, papadimitriou1998combinatorial}. Few domains illustrate as clearly how a simple model can simultaneously enable efficient computation and impose intrinsic computational limits.

The appeal of network optimization, however, extends far beyond theory. These models govern routing in transportation systems, packet forwarding in communication networks, logistics in supply chains, and the operation of energy and infrastructure networks. In each of these settings, decisions are made on a network whose link weights encode travel times, delays, costs, or capacities. The graph abstraction is compelling precisely because it mirrors the physical structure of the system itself. At the same time, the very quantities encoded by the weights are often the ones that fluctuate most.

Classical theory studies these problems under the assumption that link weights are known and fixed. In practice, this assumption is routinely violated, sometimes mildly and sometimes catastrophically. Travel times fluctuate with congestion, local disturbances spill over to neighboring streets, and disruptions propagate through interconnected links. The topology of the network remains, but its weights evolve in response to interacting flows. Consider emergency response in an urban transportation network. An ambulance must travel from a station to a hospital as quickly as possible, using a snapshot of estimated travel times inferred from sensors and historical data. If these values were static, the task would reduce to computing a shortest path. Yet the operational risk lies precisely in the fact that traffic conditions are neither static nor independent across streets. When a primary artery becomes congested, vehicles divert to adjacent roads; queues spill back across intersections; local disturbances propagate along neighboring links. A street that appeared uncongested minutes earlier may become severely delayed because flow has shifted from a nearby segment. In such systems, uncertainty is not isolated noise attached independently to edges. It is a phenomenon that moves through the network.

Robust optimization provides a principled framework for hedging against adverse realizations of uncertain parameters by replacing a nominal problem with a minimax formulation over a prescribed uncertainty set. Over the past three decades, robust optimization has developed into a powerful theoretical framework with strong duality results and tractable reformulations in many convex settings \cite{ben1998robust, ben1999robust, ben2009robust, bertsimas2011theory, gabrel2014recent}. The prevailing modeling paradigm constrains deviations of the weight vector within geometric sets in Euclidean space, such as boxes or norm balls. These sets control the magnitude and global correlation of perturbations and have enabled a broad range of computationally efficient robust counterparts. The difficulty is that computational convenience often comes with structural assumptions, and in networks those assumptions can quietly conflict with the physics of propagation. There is, in particular, a structural mismatch in many standard robust-optimization models for networks. These uncertainty sets abstract away the network’s topology. They allow perturbations to be redistributed freely across links without regard to adjacency, as if disturbances could jump from one part of the graph to another in a single step. In other words, they often treat the graph as if it were complete, even when the system is not.

Several uncertainty models have been studied specifically for network optimization \cite{aissi2009min, kouvelis2013robust,kasperski2016robust}. Interval uncertainty allows edge weights to vary within prescribed bounds, often with a budget on how many weights may deviate adversely from their nominal values \cite{bertsimas2003robust}; this leads to tractable robust formulations but ignores interaction across edges. Scenario-based models specify a finite collection of network realizations and optimize against the worst-case scenario, offering flexibility at the price of reducing an effectively continuous uncertainty to finitely many configurations \cite{kouvelis2013robust}. In many flow-driven systems, the space of plausible evolutions is effectively infinite or exponentially large, and finite scenario sets inevitably omit configurations that arise in practice. When uncertainty is allowed to range over rich infinite families, robust network optimization problems are computationally intractable in general \cite[Theorem~7]{kouvelis2013robust}. Thus we face a familiar dilemma. Models that are computationally benign are often structurally blind, and models that are structurally faithful are often algorithmically unforgiving. The natural question, then, is whether one can impose topology on uncertainty without paying for it with intractability. 

This paper argues that the dilemma is not inevitable. We propose a topology-aware uncertainty model in which perturbations are not arbitrary vectors but diffusive reallocations constrained by the directed structure of the graph. Deviations in link weights must propagate along adjacent edges and satisfy conservation constraints at nodes. The adversary can redistribute weight, but only through feasible diffusions that respect local conservation and prescribed budget constraints. This captures, at the level of the uncertainty set itself, the basic physical intuition that congestion and disruption do not teleport across a network. Rather, they flow.

Once uncertainty is forced to respect topology, the robust formulation changes structurally, and with it the relevant algorithmic questions. In convex network problems, the adversarial diffusion layer can be handled through strong duality, yielding tractable reformulations and preserving polynomial-time solvability. In combinatorial problems, however, the interaction between selecting a discrete structure and confronting a diffusive adversary can reshape the computational landscape. Seemingly minor modeling choices, such as whether diffusion is short-term or long-term and whether the adversarial budget is imposed locally or globally, can determine whether the robust problem remains tractable or becomes \textsf{NP}-hard. This leads to a structural question at the interface of robust optimization and the computational complexity of graph optimization problems:

\begin{center}
\emph{When uncertainty is constrained to diffuse along the network, when does robustification preserve algorithmic structure, and when does it create new sources of computational complexity?}
\end{center}

The answer, as we show, is not uniform. Under diffusion-based uncertainty, convex flow problems retain tractability through strong duality, while combinatorial problems exhibit a sharp and informative complexity landscape whose transitions are driven by the structure of the uncertainty set itself. We illustrate this landscape through two complementary canonical problems: shortest path, which is polynomial-time solvable in the nominal setting, and the traveling salesman problem, which is already \textsf{NP}-hard. For shortest path, short-term diffusion preserves tractability through a reduction to ordinary shortest path, whereas long-term diffusion pushes the problem across the boundary into intractability. For TSP, three of the four robust counterparts collapse back to ordinary TSP through explicit fixed-tour formulas. By formalizing diffusion-based uncertainty on directed graphs and analyzing these robust counterparts, we reveal how imposing a physically meaningful constraint on the adversary can either preserve classical algorithmic structure or create new sources of computational hardness.


\subsection{Problem formulation}
\label{subsec:prob:form}

Let $G=( V, E)$ be a directed graph with $n\coloneqq| V|$ vertices and $m\coloneqq| E|$ directed edges. We fix an arbitrary ordering of the edges so that vectors in $\mathbb{R}^m$ are indexed by $e\in  E$. Let $w\in\mathbb{R}^m_{\ge 0}$ denote the nominal edge-weight vector, where $w_e$ represents the baseline cost, delay, capacity, or other quantity of interest attached to edge $e$. The graph is assumed fixed; uncertainty enters only through structured perturbations of these edge weights. For each edge $e\in E$, we introduce two nonnegative variables,
\[
    \Delta^+_e \ge 0
    \quad\text{and}\quad
    \Delta^-_e \ge 0,
\]
which represent the amount of perturbation mass entering and leaving edge $e$, respectively. The post-diffusion edge-weight vector is defined as $w+\Delta^+-\Delta^-$, with $e^\text{th}$ component $w_e+\Delta^+_e-\Delta^-_e$ for each $e\in E$. Thus, the adversary may increase the weight of an edge by injecting perturbation mass into it and decrease it by draining mass from it. The essential modeling constraint is that perturbation mass cannot be reassigned arbitrarily across edges but must move through the directed topology of the graph. To formalize this topological structure, for each vertex $u\in  V$ let
\[
     E_{\mathrm{in}}(u) \coloneqq \{ e\in  E : e \text{ enters } u \},
    \qquad
     E_{\mathrm{out}}(u) \coloneqq \{ e\in  E : e \text{ leaves } u \}.
\]

\medskip
\noindent\textbf{Topological conservation.}
We require the diffusive uncertainty to satisfy a conservation law at every vertex,
\begin{equation}
\label{eq:conservation-vertex}
    \sum_{e\in  E_{\mathrm{in}}(u)} \Delta^-_e
    \;=\;
    \sum_{e\in  E_{\mathrm{out}}(u)} \Delta^+_e,
    \qquad \forall u\in  V.
\end{equation}
Equation \eqref{eq:conservation-vertex} enforces that perturbation mass removed from incoming edges at $u$ must reappear on outgoing edges. In transportation terms, congestion may spill over to adjacent road segments; in our model this effect is represented abstractly by node-wise conservation of perturbation mass. It is convenient to rewrite \eqref{eq:conservation-vertex} in matrix form. Define $M^+,M^-\in\mathbb{R}^{n\times m}$ by, for all $u\in V$ and $e\in E$,
\[
    (M^+)_{u,e} \coloneqq \mathbf{1}\{ e\in  E_{\mathrm{out}}(u)\},
    \qquad
    (M^-)_{u,e} \coloneqq \mathbf{1}\{ e\in  E_{\mathrm{in}}(u)\},
\]
so that conservation is equivalently expressed as
\begin{equation}
\label{eq:conservation-matrix}
    M^+\Delta^+ \;=\; M^-\Delta^-.
\end{equation}

\medskip
\noindent\textbf{Propagation regimes.}
In addition to conservation, we must specify how much perturbation mass may leave an edge. We distinguish two diffusion regimes. In a short-term (one-step) diffusion model, the outflow on each edge is limited by its baseline weight,
\begin{equation}
\label{eq:short-term}
    \Delta^- \le w,
\end{equation}
with the inequality interpreted entrywise. This captures a single-step spillover effect: perturbation mass may be redistributed to outgoing edges through a node, but an edge cannot forward more mass than it initially carries. In particular, mass that arrives at an edge cannot be forwarded again, so diffusion does not compound along longer directed paths. In a long-term (multi-step) diffusion model, outflow may be fueled by inflow,
\begin{equation}
\label{eq:long-term}
    \Delta^- \le w + \Delta^+.
\end{equation}
Here perturbation mass that arrives from upstream may accumulate and be passed further downstream. As a result, diffusion can propagate across multiple successive edges, allowing disturbances to compound along directed paths.

\medskip
\noindent\textbf{Uncertainty budgets.}
We further control the magnitude of diffusion through a budget parameter $\varepsilon\ge 0$. Let $\|\cdot\|_1$ and $\|\cdot\|_\infty$ denote the standard vector norms on $\mathbb{R}^{2m}$. Under a local ($\ell_\infty$) budget constraint, each component of $(\Delta^+,\Delta^-)$ is bounded,
\[
    \|(\Delta^+,\Delta^-)\|_\infty \le \varepsilon,
\]
which imposes a uniform per-edge cap on the amount of perturbation mass that can enter or leave any single edge. Under a global ($\ell_1$) budget constraint, the total perturbation mass is bounded,
\[
    \|(\Delta^+,\Delta^-)\|_1 \le \varepsilon,
\]
so that the uncertainty is controlled only in aggregate: the same total budget can be concentrated on a few edges or spread across many, but cannot increase overall (e.g., a limited amount of disruption that can be allocated across the network). Therefore, under an $\ell_1$ budget, conservation \eqref{eq:conservation-matrix} implies that $\sum_{e\in  E}\Delta^+_e=\sum_{e\in  E}\Delta^-_e\le \varepsilon/2$.

\medskip
\noindent\textbf{Diffusive uncertainty sets.}
Combining conservation, the diffusion regimes, and the uncertainty budgets yields four polyhedral uncertainty sets, where $\mathrm{S}$ denotes the short-term (one-step) regime and $\mathrm{L}$ the long-term (multi-step) regime:
\[
\begin{aligned}
\mathcal D^{\mathrm{S},\infty}(\varepsilon)
&\coloneqq \Bigl\{
(\Delta^+,\Delta^-)\in\mathbb{R}^{2m}_{\ge 0} :
M^+\Delta^+ = M^-\Delta^-,
\ \Delta^- \le w,
\ \|(\Delta^+,\Delta^-)\|_\infty \le \varepsilon
\Bigr\},\\
\mathcal D^{\mathrm{S},1}(\varepsilon)
&\coloneqq \Bigl\{
(\Delta^+,\Delta^-)\in\mathbb{R}^{2m}_{\ge 0} :
M^+\Delta^+ = M^-\Delta^-,
\ \Delta^- \le w,
\ \|(\Delta^+,\Delta^-)\|_1 \le \varepsilon
\Bigr\},\\
\mathcal D^{\mathrm{L},\infty}(\varepsilon)
&\coloneqq \Bigl\{
(\Delta^+,\Delta^-)\in\mathbb{R}^{2m}_{\ge 0} :
M^+\Delta^+ = M^-\Delta^-,
\ \Delta^- \le w+\Delta^+,
\ \|(\Delta^+,\Delta^-)\|_\infty \le \varepsilon
\Bigr\},\\
\mathcal D^{\mathrm{L},1}(\varepsilon)
&\coloneqq \Bigl\{
(\Delta^+,\Delta^-)\in\mathbb{R}^{2m}_{\ge 0} :
M^+\Delta^+ = M^-\Delta^-,
\ \Delta^- \le w+\Delta^+,
\ \|(\Delta^+,\Delta^-)\|_1 \le \varepsilon
\Bigr\}.
\end{aligned}
\]
Each uncertainty set is a bounded polyhedron and therefore convex and compact, a structural property that will play a central role in the analysis of convex problems.

\medskip
\noindent\textbf{Diffusion-robust optimization.} 
We study optimization problems on $G$ in which a decision variable $f$ belongs to a feasible set $\mathcal F$ encoding a network structure, such as flows, paths, or tours, while an adversary selects a diffusion $(\Delta^+,\Delta^-)$ from one of the uncertainty sets defined above. This leads to the diffusion-robust minimax formulation
\begin{equation}
    \min_{f\in \mathcal F}\ \max_{(\Delta^+,\Delta^-)\in \mathcal D(\varepsilon)}\ 
    g\left(f,\Delta^+,\Delta^-\right),
    \label{eq:minimax-general}
\end{equation}
where $\mathcal D(\varepsilon)$ denotes one of the four regimes. In additive cost settings,
\[
    g(f,\Delta^+,\Delta^-)
    =
    f^\top\bigl(w+\Delta^+-\Delta^-\bigr),
\]
so that the adversary worsens performance by transporting perturbation mass through the network subject to conservation and budget constraints. For linear network-flow models, such as minimum-cost flow, standard linear programming duality can be used to dualize the inner maximization over $\mathcal D(\varepsilon)$. Thus, across all four diffusion regimes, the resulting robust counterparts remain linear programs and can be solved in polynomial time. In this paper, we instantiate \eqref{eq:minimax-general} for two canonical combinatorial graph problems, shortest path and traveling salesman, and analyze how the structure of $\mathcal D(\varepsilon)$ shapes the boundary between polynomial-time solvability and \textsf{NP}-hardness.


\subsection{Contributions}
\label{subsec:contributions}

The main contributions of the paper can be summarized as follows.

\begin{itemize}
\item[(i)] \textbf{Diffusive uncertainty model.}
We introduce a new diffusion-based uncertainty model for edge weights on directed graphs: perturbations propagate through adjacency subject to node-wise conservation, under two propagation regimes (short-term vs.\ long-term) and two budget geometries (local $\ell_\infty$ vs.\ global $\ell_1$). The resulting uncertainty sets are polyhedral and compact.

\item[(ii)] \textbf{Shortest path.}
For diffusion-robust shortest path (Diff-RSP), we establish a propagation-driven complexity transition:
\begin{itemize}
    \item[1.] \emph{Short-term diffusion.} 
    Diff-RSP is polynomial-time solvable under $\mathcal D^{\mathrm{S},\infty}(\varepsilon)$ and $\mathcal D^{\mathrm{S},1}(\varepsilon)$. In these regimes, the inner maximization admits an explicit closed form, reducing the robust problem to a constant number of ordinary shortest-path instances on precomputed edge weights.

    \item[2.] \emph{Long-term diffusion.} 
    Diff-RSP is \textsf{NP}-hard under both $\mathcal D^{\mathrm{L},\infty}(\varepsilon)$ and $\mathcal D^{\mathrm{L},1}(\varepsilon)$. In both cases, multi-step propagation allows perturbation mass to reach the chosen path through the surrounding network, fundamentally altering tractability. Under the global-budget regime $\mathcal D^{\mathrm{L},1}(\varepsilon)$, this propagation has a transportation structure that encodes exposure to the surrounding network, yielding a polynomial-time reduction from \emph{Most Secluded Path}. Under the local-budget regime $\mathcal D^{\mathrm{L},\infty}(\varepsilon)$, it instead creates a binary activation mechanism, yielding a polynomial-time reduction from \emph{Minimum Satisfiability (MinSAT)}.

\end{itemize}

\item[(iii)] \textbf{Traveling salesman.}
For diffusion-robust traveling salesman (Diff-RTSP), we show that robustification often adds no computational complexity beyond ordinary TSP, but that one regime exhibits a more intricate fixed-tour evaluation problem:
\begin{itemize}
    \item[1.] \emph{Exact reductions to ordinary TSP.}
    Under $\mathcal D^{\mathrm{S},\infty}(\varepsilon)$, $\mathcal D^{\mathrm{S},1}(\varepsilon)$, and $\mathcal D^{\mathrm{L},1}(\varepsilon)$, Diff-RTSP is polynomial-time equivalent to ordinary TSP. The fixed-tour worst-case value collapses either to the cost of the same tour under precomputed worst-case edge weights or to a scalar capped expression depending only on the tour's nominal weight.

    \item[2.] \emph{The long-term local-budget regime.}
    Under $\mathcal D^{\mathrm{L},\infty}(\varepsilon)$, multi-step propagation and local edgewise caps interact, so the fixed-tour adversarial problem does not generally collapse to an ordinary TSP objective. We show that a natural upper bound obtained from two ordinary TSP instances can be strict, even on a complete directed graph, and we bracket the optimal robust value between quantities computable from ordinary TSP instances.
\end{itemize}
\end{itemize}


\subsection{Related work}
\label{subsec:related:work}

Robust graph optimization is generally studied under two optimality criteria. The most common is the \emph{min--max criterion}, which minimizes the worst realized cost of the chosen solution over the uncertainty set. This is the criterion adopted in this paper, as reflected in the formulation \eqref{eq:minimax-general}. A second widely studied criterion is \emph{min--max regret}, where the objective is the worst excess cost relative to the solution that would have been optimal in hindsight. We refer to \cite{aissi2009min,kouvelis2013robust,kasperski2016robust} for surveys of robust combinatorial optimization, including min--max regret formulations.

Within the min--max setting, the models most closely related to ours are \emph{static, single-stage} robust formulations, in which a feasible solution is chosen before uncertainty is realized and is evaluated against the worst admissible realization. Existing models in this class can be organized by how they specify the uncertainty set. A large part of the robust graph-optimization literature is built around \emph{scenario}, \emph{interval}, and \emph{budgeted} uncertainty. Scenario models represent uncertainty by a finite collection of possible realizations; interval models allow each coefficient to vary within a prescribed range; and, among interval-based models, budgeted models restrict how many edge weights can deviate simultaneously from their nominal values. Recent work has also considered alternative descriptions of uncertainty, including \emph{ellipsoidal}, \emph{data-driven}, and \emph{distributionally robust} uncertainty sets, which impose geometric or statistical structure on the set of admissible edge-weight vectors. These approaches control the size, dependence, or statistical plausibility of the uncertain weight vector. The model studied here is also static and single-stage, but addresses a different structural feature: \emph{uncertainty that propagates through the topology of the graph}. Rather than treating uncertainty as a topology-agnostic set of admissible edge-weight vectors, our model requires perturbations to arise from feasible diffusions of perturbation mass through adjacent edges, subject to node-wise conservation.

\medskip
\noindent\textbf{Robust Shortest Path.}
Among robust graph-optimization problems, robust shortest path is one of the canonical models and has received sustained attention. The literature is now broad enough that several reference points are available: \cite{filippi2025robust} surveys robust and distributionally robust shortest-path formulations, \cite[Chapter~7]{goerigk2024introduction} summarizes complexity and approximability results for robust shortest path under different criteria and uncertainty sets, and \cite{goerigk2024benchmarking} develops benchmark instances for comparing robust discrete-optimization models, including robust shortest-path variants. We therefore give only a selection of representative examples. The most common uncertainty models in this literature are scenario, interval, and budgeted uncertainty. Scenario-based robust shortest-path models and exact methods are studied in \cite{murthy1992solving, yu1998robust, bruni2010enhanced, xing2013reformulation, duque2019exact}. For the plain min--max criterion, independent interval uncertainty reduces to an ordinary shortest-path problem with edge weights set to their upper bounds; interval uncertainty becomes substantially richer under min--max regret and relative-robustness formulations, leading to mixed-integer formulations \cite{karasan2001robust}, general interval-data complexity results \cite{averbakh2004interval}, exact and branch-and-bound algorithms \cite{montemanni2004exact, montemanni2004branch, montemanni2005robust}, computational studies \cite{zielinski2004computational}, and reduction or special-graph results \cite{kasperski2006robust, catanzaro2011reduction}. Budgeted uncertainty is studied in \cite{bertsimas2003robust}, while robust shortest-path models with additional resource or feasibility constraints are studied in \cite{kwon2013robust, dipuglia2019resource}. Other uncertainty descriptions include ellipsoidal uncertainty \cite{alves2015robust,goldberg2025smooth}, data-driven robust shortest path \cite{shahabi2015robust}, and distributionally robust shortest path \cite{cheng2013distributionally, cheng2016new, wang2020wasserstein, ketkov2021approach, ketkov2023multistage}. Additional robust shortest-path variants include multiobjective formulations \cite{chassein2019algorithms}, algorithmic and computational developments \cite{raith2018extensions, bertsekas2019robust, hansknecht2018fast}, recourse and recoverability \cite{golovin2015improved, li2015large, busing2012recoverable, jackiewicz2025computational}, and dynamic robust shortest path \cite{xu2020dynamic}. 

\medskip
\noindent\textbf{Robust Traveling Salesman (and Vehicle Routing).}
Although smaller than the robust shortest-path literature, this line of work follows a similar modeling pattern: uncertainty is typically imposed on edge costs, travel times, or the ability to revise a tour after uncertainty is observed. We again give only a selection of representative examples. Interval-data robust TSP has been studied through theoretical properties, formulations, and exact and heuristic algorithms \cite{montemanni2007robust}. Recoverable robust variants, in which the tour may be modified after uncertainty is observed, are studied in \cite{chassein2016recoverable, goerigk2021recoverable}. Other formulations include Wasserstein distributionally robust Euclidean TSP \cite{carlsson2018wasserstein} and robust TSP with time windows under knapsack-constrained travel-time uncertainty \cite{bartolini2021robust}. Robust-regret algorithms have also been developed for \textsf{NP}-hard graph optimization problems, including TSP and Steiner tree \cite{ganesh2023robust}. The closely related vehicle-routing literature studies robust variants under richer operational constraints, including uncertain demands \cite{ordonez2010robust, gounaris2013robust}, travel times and service times \cite{hu2018robust, munari2019robust}, and distributional ambiguity \cite{ghosal2020distributionally}. These works show that robust routing problems depend strongly on which operational quantities are uncertain and whether decisions can be revised. Our TSP results add a different axis: even when the underlying combinatorial problem is already \textsf{NP}-hard, the topology imposed on the uncertainty set determines whether the diffusion-robust counterpart admits an exact reduction to ordinary TSP or instead requires a separate analysis of the fixed-tour adversarial problem.

\medskip
\noindent\textbf{Other Robust Graph and Network Models.}
Robust graph optimization also includes models in which uncertainty affects demands, capacities, recourse decisions, or the availability of network components. Two-stage robust network flow and design are studied in \cite{atamturk2007two, nasrabadi2013robust}, incremental and recoverable robustness in network problems is studied in \cite{cseref2009incremental, liebchen2009concept}, and robust or adaptive flow models with node or arc failures are studied in \cite{bertsimas2013robust}. These works change the source of uncertainty, the timing of decisions, or the information available after uncertainty is revealed. The present paper is complementary: our model remains static and single-stage, but makes the directed topology of the graph part of the uncertainty set itself through diffusion and node-wise conservation.


\subsection{Organization and notation}
\label{subsec:organization:notation}

\noindent\textbf{Organization.}
Section~\ref{sec:shortest-path} studies Diff-RSP, with Section~\ref{subsec:rsp-poly} giving algorithms for the short-term regimes and Sections~\ref{subsec:rsp-np-local} and~\ref{subsec:rsp-np-global} proving \textsf{NP}-hardness for the long-term local- and global-budget regimes via reductions from MinSAT and Most Secluded Path, respectively. Section~\ref{sec:tsp} studies Diff-RTSP, with Sections~\ref{subsec:tsp-short-infty},~\ref{subsec:tsp-global}, and~\ref{subsec:tsp-long-infty} covering the short-term local-budget regime, the two global-budget regimes, and the long-term local-budget regime, respectively. All proofs are collected in the appendix.

\medskip
\noindent\textbf{Notation.}
For a vector $x\in\mathbb{R}^m$, we write $\|x\|_1$ and $\|x\|_\infty$ for the standard $\ell_1$ and $\ell_\infty$ norms. Throughout, $\le_p$ denotes polynomial-time reducibility, and $\textsf{P}$ and $\textsf{NP}$ denote the standard complexity classes.


\section{Diffusion-Robust Shortest Path}
\label{sec:shortest-path}

Shortest path is a canonical problem in network optimization: it is polynomial-time solvable, admits several classical combinatorial algorithms (e.g., Dijkstra, Bellman--Ford, and DAG shortest-path algorithms), and can be formulated as a minimum-cost flow problem with a tight linear relaxation. This makes it a natural benchmark for understanding how diffusion-based uncertainty reshapes computational structure. Fix source and sink nodes $s,t\in V$. Let $B\in\mathbb{R}^{n\times m}$ denote the node--edge incidence matrix of $G$, and let $b^{s,t}\coloneqq e_s-e_t\in\mathbb{R}^n$, where $e_u$ is the $u^\text{th}$ standard basis vector. We define
\[
    \mathcal P_{s,t}\coloneqq\{f\in\{0,1\}^m : Bf=b^{s,t}\}.
\]
Although a vector $f\in\mathcal P_{s,t}$ may contain directed cycles in addition to an $s$--$t$ path, the diffusion constraints imply that every realized edge weight $w_e+\Delta^+_e-\Delta^-_e$ is nonnegative under each regime. Thus deleting any directed cycle preserves feasibility and weakly decreases the realized cost for every feasible perturbation. Consequently, the robust optimum over $\mathcal P_{s,t}$ is attained by a cycle-free $s$--$t$ path, so we may identify feasible solutions with directed $s$--$t$ paths throughout. The diffusion-robust shortest path problem, abbreviated Diff-RSP, is
\[
    \mathrm{OPT}_{\mathrm{RSP}}(w,\mathcal D(\varepsilon))
    :=
    \min_{f\in\mathcal P_{s,t}}\; \max_{(\Delta^+,\Delta^-)\in \mathcal D(\varepsilon)}
    f^\top\bigl(w+\Delta^+-\Delta^-\bigr),
\]
where $\mathcal D(\varepsilon)$ denotes one of the four diffusion uncertainty sets. Unlike the nominal case, robustification introduces an adversarial value function that is generally not edge-separable, so the objective is no longer a linear path cost. In particular, the worst-case perturbation can couple edges through the diffusion constraints, and the standard minimum-cost-flow formulation no longer applies directly. As the next theorem shows, this coupling creates a sharp complexity split: short-term diffusion remains polynomial-time solvable, whereas long-term diffusion becomes \textsf{NP}-hard.

\begin{theorem}[Complexity of Diff-RSP]
\label{thm:diff-rsp}
Diff-RSP exhibits a complexity transition driven by propagation depth:
\begin{enumerate}[label=(\roman*)]
    \item Under $\mathcal D^{\mathrm{S},\infty}(\varepsilon)$, Diff-RSP is polynomial-time solvable.
    \item Under $\mathcal D^{\mathrm{S},1}(\varepsilon)$, Diff-RSP is polynomial-time solvable.
    \item Under $\mathcal D^{\mathrm{L},\infty}(\varepsilon)$, Diff-RSP is \textsf{NP}-hard.
    \item Under $\mathcal D^{\mathrm{L},1}(\varepsilon)$, Diff-RSP is \textsf{NP}-hard.
\end{enumerate}
\end{theorem}

The proof proceeds by separating the short-term and long-term regimes. In the short-term cases, the inner adversarial problem admits an explicit evaluation, which yields exact reductions to ordinary shortest-path instances with modified edge costs. In the long-term cases, multi-step diffusion creates path-level coupling and leads to \textsf{NP}-hardness through two different mechanisms: the global-budget regime encodes exposure to the surrounding network via \emph{Most Secluded Path}, whereas the local-budget regime encodes clause activation via \emph{Minimum Satisfiability} (MinSAT).


\subsection{Algorithms for the polynomial cases}
\label{subsec:rsp-poly}

The two short-term regimes are tractable because the short-term diffusion constraint $\Delta^- \le w$ ensures that the amount leaving any edge is controlled solely by its nominal weight, rather than by perturbation mass that has newly arrived there. Hence diffusion is intrinsically one-step: perturbation mass may be removed from edges entering a node and reassigned to edges leaving that node, but once it has been reassigned, it cannot continue propagating further downstream. In particular, the adversary cannot create recursive multi-step amplification along an $s$--$t$ path. For any fixed path, the inner maximization therefore reduces to a local one-step redistribution effect that can be evaluated exactly.

The algorithm exploits this one-step structure through a preprocessing step that is shared by both polynomial cases. The central idea is to encode, in advance, the worst local diffusion that can occur when a path passes through a node. Consider a node $u$ traversed by a candidate $s$--$t$ path, and let $e_u^{\mathrm{in}}$ and $e_u^{\mathrm{out}}$ denote the edges on which the path enters and leaves $u$, respectively. Once the path reaches $u$ through $e_u^{\mathrm{in}}$, the incoming neighborhood of $u$ determines the maximum one-step perturbation that can be transferred through $u$ and used to worsen the next step. Thus, the worst local effect of passing through $u$ is determined entirely by the edges entering $u$.

\begin{figure}[t]
    \centering
    \includegraphics[width=0.55\linewidth]{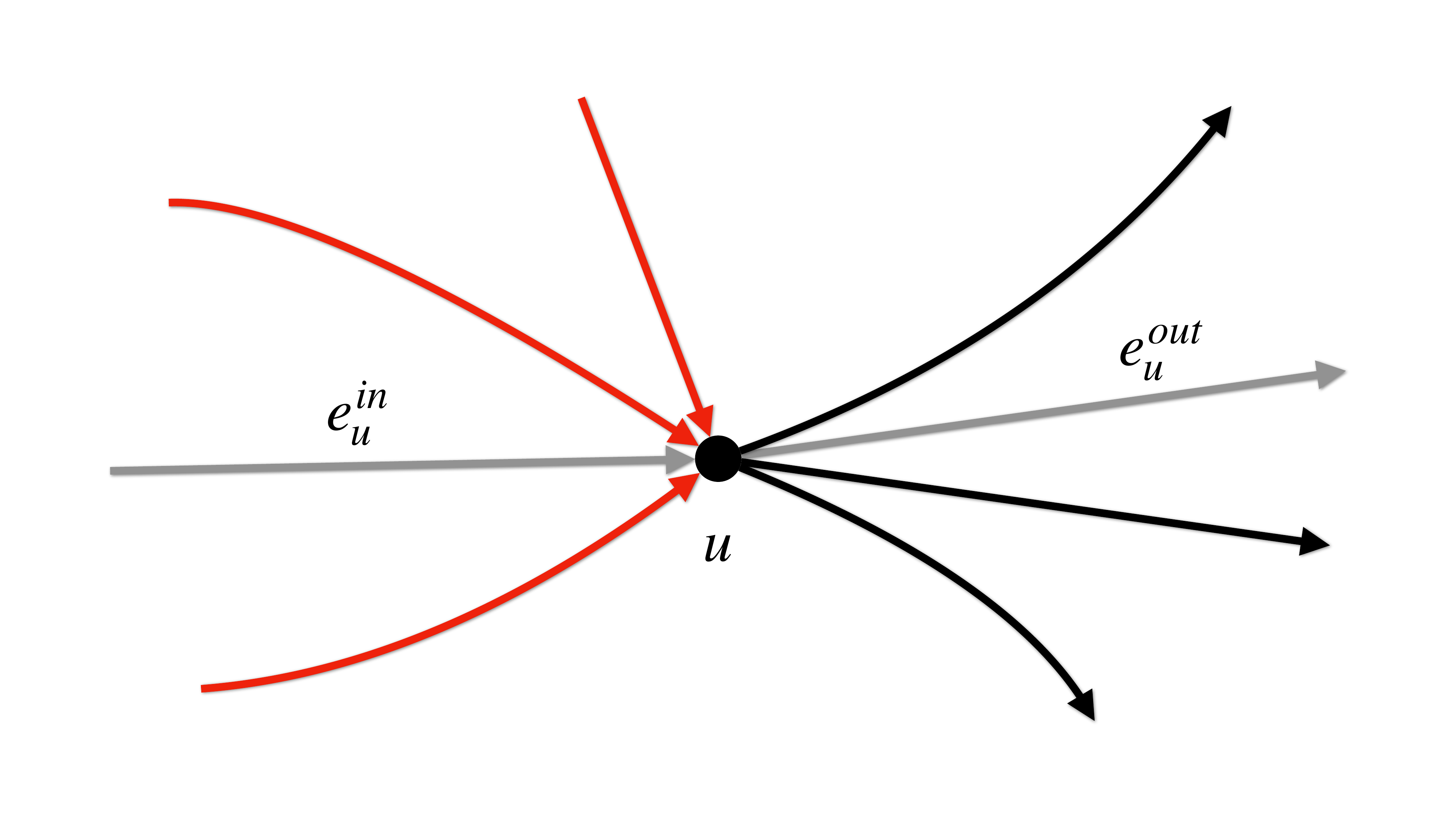}
    \caption{Local diffusion effect at a path node $u$. The gray edges $e_u^{\mathrm{in}}$ and $e_u^{\mathrm{out}}$ are the edges used by the path to enter and leave $u$. The red edges indicate other incoming edges of $u$. Under short-term diffusion, perturbation mass available at the incoming edges of $u$ may be transferred through $u$ and used to worsen the next step $e_u^{\mathrm{out}}$. However, the amount that can leave any edge is bounded solely by its nominal weight, so downstream propagation cannot be amplified by newly arrived perturbation mass.}
    \label{fig:rsp-local-diffusion}
\end{figure}

Figure~\ref{fig:rsp-local-diffusion} illustrates this mechanism. The physical effect of diffusion is realized on the outgoing path edge $e_u^{\mathrm{out}}$: after the path arrives at $u$, the adversary can use perturbation mass available at $u$ to make the next move more expensive. The key algorithmic step is to account for this future increase \emph{in advance}. Instead of attaching the corresponding surcharge directly to $e_u^{\mathrm{out}}$, we record an equivalent surcharge on the current edge $e_u^{\mathrm{in}}$. This accounting step is the key device that makes the algorithm possible. For every $s$--$t$ path, it produces the same total worst-case cost as the original diffusion model, but it assigns each local diffusion penalty to an edge already chosen by the path. Thus the cost of a path can be computed by simply summing precomputed edge weights, rather than tracking how perturbation mass becomes available at successive vertices. In this way, the robust instance is converted into a single precomputed worst-case graph with modified edge weights.

In both polynomial cases, we first compute a path-independent bound on the one-step perturbation that can be transferred through each node, and then use these node-level quantities to define the surcharges in the worst-case graph. Once these surcharges have been computed, the worst-case value of any fixed $s$--$t$ path can be evaluated by summing modified edge costs, together with an additive term that accounts for the diffusion contribution at the source node $s$. Consequently, Diff-RSP reduces to ordinary shortest-path computations. Under the local $\ell_\infty$ budget, this yields a single shortest-path instance with modified edge weights. Under the global $\ell_1$ budget, the same preprocessing applies, but the aggregate budget constraint produces two competing separable path costs, so the robust value is obtained from two shortest-path computations.

We now formalize the preprocessing step. For each node $u\in V$, define the one-step transfer bound
\[
    T_u \;\coloneqq\; \sum_{e\in  E_{\mathrm{in}}(u)} \min\{\varepsilon,w_e\}.
\]
This is the maximum perturbation mass that can be collected from the incoming neighborhood of $u$ under short-term diffusion. Next, for each edge $e=(v,u)\in E$, define the surcharge
\[
    \chi_e
    \;\coloneqq\;
    \begin{cases}
        0, & u=t,\\[0.3em]
        \min\bigl\{\varepsilon,\; T_u-\min\{\varepsilon,w_e\}\bigr\} \;=\; \min\bigl\{\varepsilon,\; \sum_{e^\prime \in  E_{\mathrm{in}}(u)\setminus\{e\}} \min\{\varepsilon,w_{e^\prime}\}\bigr\}, & u\neq t.
    \end{cases}
\]
Thus, $\chi_e$ is the surcharge assigned to edge $e$ in order to account in advance for the worst one-step diffusion that may occur when the path subsequently passes through the head node $u$ of $e$. We also define the source correction term
\[
    c_s \;\coloneqq\; \min\{\varepsilon,T_s\},
\]
which accounts for the diffusion contribution generated at the source node $s$. Finally, let
\[
    w_e^{\mathrm{wc}} \;\coloneqq\; w_e+\chi_e,
    \qquad e\in E,
\]
and write $w^{\mathrm{wc}}\in\mathbb{R}^m_{\ge 0}$ for the resulting edge-weight vector of the precomputed \emph{worst-case graph}.

\begin{algorithm}[t]
\caption{Diff-RSP under short-term diffusion}
\label{alg:diff-rsp-short-term}
\begin{algorithmic}[1]
\Require A directed graph $G=( V, E)$ with edge weights $w\in\mathbb{R}_{\ge 0}^m$, source $s\in V$, sink $t\in V$, budget $\varepsilon\ge 0$, and regime $\mathsf{Reg}\in\{\mathcal D^{\mathrm{S},\infty}(\varepsilon),\mathcal D^{\mathrm{S},1}(\varepsilon)\}$.
\Ensure An optimal path-incidence vector $f^\star$ and the corresponding optimal robust value $\mathrm{OPT}_{\mathrm{RSP}}$.

\Statex \textbf{Preprocessing}
\State Compute the node-level transfer bounds $\{T_u\}_{u\in V}$, the edge surcharges $\{\chi_e\}_{e\in E}$, the worst-case graph weights $w^{\mathrm{wc}}$, and the source correction term $c_s$ as defined above.

\Statex \textbf{Case 1: local budget $\mathcal D^{\mathrm{S},\infty}(\varepsilon)$}
\If{$\mathsf{Reg}=\mathcal D^{\mathrm{S},\infty}(\varepsilon)$}
    \State Compute a shortest $s$--$t$ path in $G$ with edge weights $w^{\mathrm{wc}}$.
    \State Let $f^\star$ be its incidence vector.
    \State Set $\mathrm{OPT}_{\mathrm{RSP}}\gets (w^{\mathrm{wc}})^\top f^\star+c_s$.
    \State \Return $(f^\star,\mathrm{OPT}_{\mathrm{RSP}})$.
\EndIf

\Statex \textbf{Case 2: global budget $\mathcal D^{\mathrm{S},1}(\varepsilon)$}
\If{$\mathsf{Reg}=\mathcal D^{\mathrm{S},1}(\varepsilon)$}
    \State Compute a shortest $s$--$t$ path in $G$ with edge weights $w^{\mathrm{wc}}$.
    \State Let $f^{(1)}$ be its incidence vector.
    \State Set $\mathrm{VAL}^{(1)}\gets (w^{\mathrm{wc}})^\top f^{(1)}+T_s$.

    \State Compute a shortest $s$--$t$ path in $G$ with edge weights $w$.
    \State Let $f^{(0)}$ be its incidence vector.
    \State Set $\mathrm{VAL}^{(0)}\gets w^\top f^{(0)}+\varepsilon/2$.

    \If{$\mathrm{VAL}^{(1)}\le \mathrm{VAL}^{(0)}$}
        \State Set $f^\star\gets f^{(1)}$ and $\mathrm{OPT}_{\mathrm{RSP}}\gets \mathrm{VAL}^{(1)}$.
    \Else
        \State Set $f^\star\gets f^{(0)}$ and $\mathrm{OPT}_{\mathrm{RSP}}\gets \mathrm{VAL}^{(0)}$.
    \EndIf
    \State \Return $(f^\star,\mathrm{OPT}_{\mathrm{RSP}})$.
\EndIf
\end{algorithmic}
\end{algorithm}

Algorithm~\ref{alg:diff-rsp-short-term} separates the computation into a local preprocessing stage and a shortest-path stage. Assuming the graph is represented by incoming adjacency lists, the preprocessing stage computes the node-level transfer bounds and the corresponding edge surcharges that define the precomputed worst-case graph, and runs in $O(| E|+| V|)$ time: each edge contributes once to the computation of the node bound at its head, and once these node-level quantities are available, each surcharge is obtained in constant time. After preprocessing, the robust optimization problem is reduced to standard shortest-path computations. The only distinction between the two short-term regimes lies in the final selection step: under $\mathcal D^{\mathrm{S},\infty}(\varepsilon)$, the algorithm solves a single shortest-path instance on the worst-case graph with edge weights $w^{\mathrm{wc}}$, whereas under $\mathcal D^{\mathrm{S},1}(\varepsilon)$, it compares the values of two shortest-path solutions corresponding to the two candidate separable objectives.

The next proposition shows that Algorithm~\ref{alg:diff-rsp-short-term} returns an optimal solution in both short-term regimes and, in particular, establishes assertions~(i) and~(ii) of Theorem~\ref{thm:diff-rsp}.

\begin{proposition}[Optimality and complexity of Algorithm~\ref{alg:diff-rsp-short-term}]
\label{prop:diff-rsp-short-term-algorithm}
Algorithm~\ref{alg:diff-rsp-short-term} returns an optimal path-incidence vector $f^\star$ and the corresponding optimal robust value $\mathrm{OPT}_{\mathrm{RSP}}(w,\mathcal D(\varepsilon))$ for Diff-RSP under both short-term regimes. In particular,
\[
    \mathrm{OPT}_{\mathrm{RSP}}(w,\mathcal D(\varepsilon))
    =
    \min_{f\in\mathcal P_{s,t}}
    \max_{(\Delta^+,\Delta^-)\in\mathcal D(\varepsilon)}
    f^\top\bigl(w+\Delta^+-\Delta^-\bigr),
\]
when $\mathcal D(\varepsilon)$ is either $\mathcal D^{\mathrm{S},\infty}(\varepsilon)$ or $\mathcal D^{\mathrm{S},1}(\varepsilon)$. In both cases, the running time is $O(| E|+| V|\log| V|)$.
\end{proposition}

Proposition~\ref{prop:diff-rsp-short-term-algorithm} identifies the precise source of tractability in the short-term regimes: the adversarial effect can be absorbed into a path-independent reweighting of the graph. This is a genuinely structural property of one-step diffusion, not merely a consequence of small budget size. In particular, the robust objective remains compatible with the combinatorial structure of shortest path because the worst-case interaction can be encoded locally and precomputed once for all candidate paths. This path-independence will be the key feature that fails in the long-term regimes, where diffusion can propagate recursively and no analogous fixed worst-case graph exists.


\subsection{\textsf{NP}-hardness under long-term local budget}
\label{subsec:rsp-np-local}

The long-term local $\ell_\infty$-budget regime becomes hard because the adversary can \emph{move perturbation mass over multiple steps and route it onto many different edges used by the chosen path}. This is one of the mechanisms that breaks the short-term ``worst-case graph'' reduction. Under short-term diffusion, the constraint $\Delta^- \le w$ allows mass to be removed only from edges that already have nominal weight; as a result, the worst-case increase incurred when a path passes through a node remains a one-step, node-local effect that can be encoded in advance. Under long-term diffusion, the constraint becomes $\Delta^- \le w+\Delta^+$, so mass that reaches an edge can be forwarded further downstream. In the local-budget regime, this forwarding is not charged by a global transportation budget: each edge-wise perturbation is bounded, but many different diffusion routes can be active at the same time. Thus, the robust cost of a candidate path can depend on which upstream sources have diffusion routes into the particular edges selected by the path. We use this path-dependent reachability effect to encode the clauses of a Boolean formula.

The Boolean optimization problem we use is \emph{Minimum Satisfiability} (MinSAT). Informally, MinSAT asks for a truth assignment that satisfies as few clauses as possible.

\begin{definition}[Minimum Satisfiability]
\label{def:minsat}
An instance of MinSAT consists of a Boolean formula
\[
    \Phi = C_1 \wedge C_2 \wedge \cdots \wedge C_m
\]
over variables $x_1,\ldots,x_n$, where each clause $C_j$ is an OR of literals, and each literal is either a variable $x_i$ or its negation $\neg x_i$. Given a truth assignment $A$, let
\[
    \mathrm{sat}(A)
    \;\coloneqq\;
    \bigl|\{j\in\{1,\ldots,m\}: C_j \text{ is satisfied by } A\}\bigr|.
\]
The MinSAT problem asks to find a truth assignment $A$ minimizing $\mathrm{sat}(A)$.
\end{definition}

MinSAT was introduced by Kohli, Krishnamurti, and Mirchandani~\cite{kohli1994minimum}, who showed that it is \textsf{NP}-hard even when each clause contains at most two literals.

The connection to long-term local diffusion is now intuitive. The chosen path will encode a truth assignment, and the construction will introduce one perturbation source associated with each clause $C_j$. The construction is designed so that a clause source creates one additional unit contribution to the robust cost precisely when the corresponding clause is satisfied by the assignment encoded by the path. Clause sources corresponding to unsatisfied clauses do not create an additional contribution beyond the baseline. Thus, the worst-case long-term local-budget cost of the corresponding path equals a fixed baseline, independent of the assignment, plus $\mathrm{sat}(A)$. Consequently, minimizing the robust cost of the chosen path is equivalent to choosing an assignment that minimizes the number of satisfied clauses. We now give the construction.

\medskip
\noindent\textbf{Construction of $G'$.}
The construction replaces the Boolean formula $\Phi$ by a sequence of \emph{variable gadgets}. Each variable gadget offers two directed branches: a true branch $T_i$, corresponding to setting $x_i=\mathrm{true}$, and a false branch $F_i$, corresponding to setting $x_i=\mathrm{false}$. The edges on these branches represent literal occurrences in the formula. The construction also introduces two types of source edges. \emph{Clause source edges} encode the clauses of $\Phi$: the source edge associated with clause $C_j$ can route mass into the branches corresponding to the literals contained in $C_j$. \emph{Blocker source edges} create a fixed baseline contribution, one per variable gadget. This baseline is needed because a clause source corresponding to an unsatisfied clause may still travel through an unchosen branch and meet the chosen path at the next blocking edge. The blocking edges ensure that such leakage is absorbed by an assignment-independent baseline and cannot create an additional clause contribution. All source edges have nominal weight $1$, and all remaining edges have nominal weight $0$. Figure~\ref{fig:rsp-linf-construction} illustrates the variable gadgets, the blocker source edges, and the clause-source connectors. The graph $G'=(V',E')$ is constructed as follows:

\begin{figure}[t]
    \centering
    \includegraphics[width=0.95\linewidth]{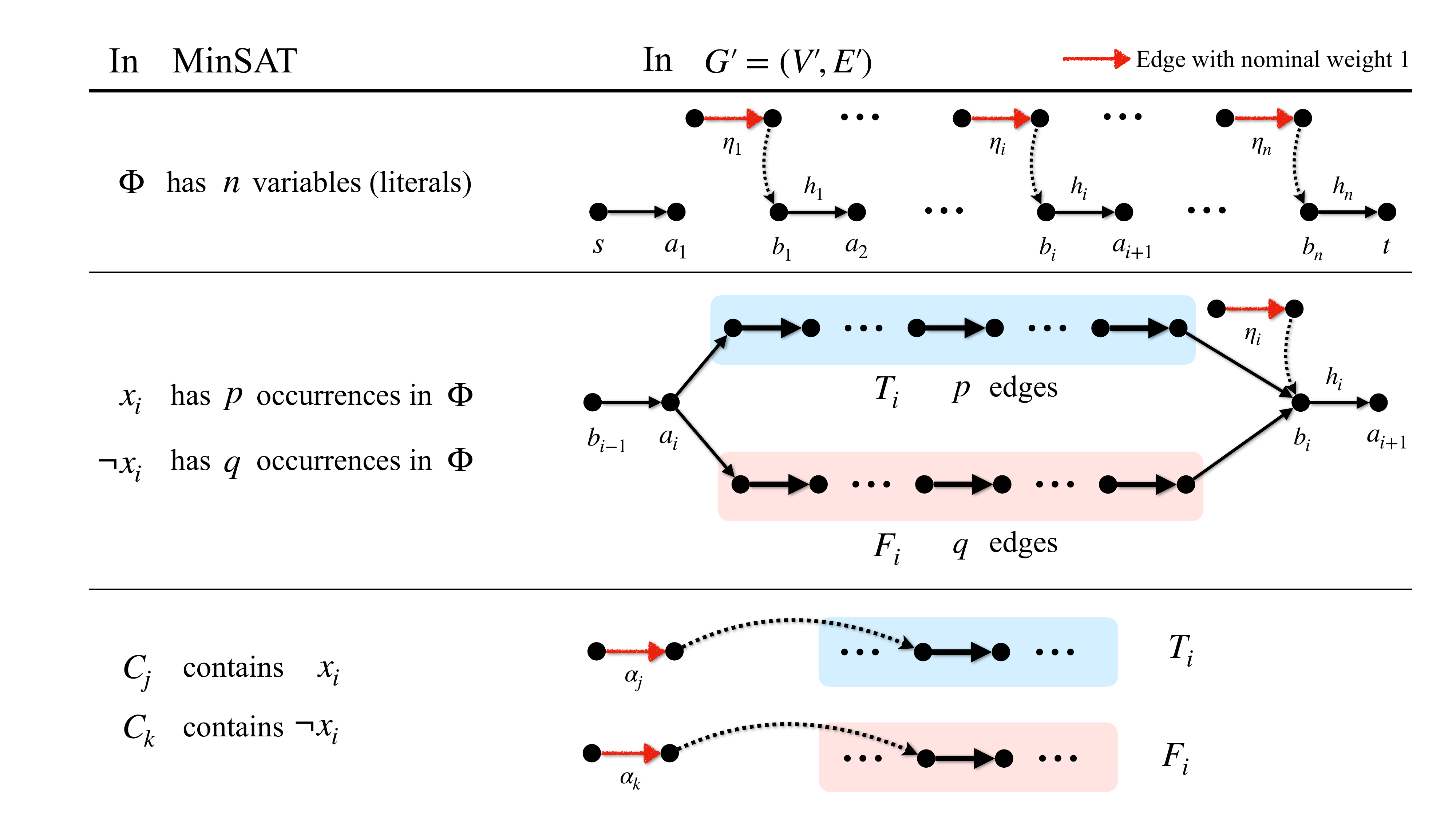}
    \caption{Construction of $G'=(V',E')$ from a MinSAT instance $\Phi=C_1\wedge\cdots\wedge C_m$. Top: the variable gadgets are arranged in series, and each gadget has a blocker source edge $\eta_i$ that can route mass to the blocking edge $h_i$. Middle: variable $x_i$ is represented by two branches from $a_i$ to $b_i$. The true branch $T_i$ contains one edge for each occurrence of $x_i$, and the false branch $F_i$ contains one edge for each occurrence of $\neg x_i$. Bottom: a clause source edge $\alpha_j$ connects to the literal-occurrence edges corresponding to the literals contained in $C_j$. Edges drawn in red have nominal weight $1$; all other edges have nominal weight $0$.}
    \label{fig:rsp-linf-construction}
\end{figure}

\begin{itemize}
    \item \textbf{Variable backbone.}
    Create vertices $s,t$ and, for each variable $x_i$, two vertices $a_i$ and $b_i$. Add a zero-weight edge from $s$ to $a_1$. For each $i=1,\ldots,n-1$, add a zero-weight edge
    \[
        h_i \coloneqq (b_i,a_{i+1}),
    \]
    and add the final zero-weight edge
    \[
        h_n \coloneqq (b_n,t).
    \]
    We call $h_1,\ldots,h_n$ the \emph{blocking edges}. Every directed $s$--$t$ path contains all blocking edges $h_1,\ldots,h_n$, one after each variable gadget.

    \item \textbf{Literal branches.}
    For each variable $x_i$, create two internally node-disjoint directed paths from $a_i$ to $b_i$. The first path is denoted by $T_i$ and corresponds to setting $x_i=\mathrm{true}$; the second path is denoted by $F_i$ and corresponds to setting $x_i=\mathrm{false}$. If $x_i$ has $p_i$ occurrences in $\Phi$, then $T_i$ contains $p_i$ distinguished \emph{literal-occurrence edges}, one for each occurrence of $x_i$. If $\neg x_i$ has $q_i$ occurrences in $\Phi$, then $F_i$ contains $q_i$ distinguished literal-occurrence edges, one for each occurrence of $\neg x_i$. These occurrence edges are private to their branches, and their tails are internal branch vertices, not the common endpoints $a_i$ or $b_i$. If $x_i$ has no occurrences, then $T_i$ is represented by a single zero-weight auxiliary edge with no clause connector; similarly, if $\neg x_i$ has no occurrences, then $F_i$ is represented by a single zero-weight auxiliary edge with no clause connector. All edges in $T_i$ and $F_i$ have nominal weight $0$.

    \item \textbf{Clause source edges.}
    For each clause $C_j$, create a dummy edge
    \[
        \alpha_j \coloneqq (d_j,c_j),
    \]
    with nominal weight $1$. For each occurrence of a literal in $C_j$, add a zero-weight connector from $c_j$ to the tail of the corresponding literal-occurrence edge inside the appropriate branch. Thus, if $C_j$ contains an occurrence of $x_i$, then this connector enters the corresponding occurrence edge on the true branch $T_i$; if $C_j$ contains an occurrence of $\neg x_i$, then this connector enters the corresponding occurrence edge on the false branch $F_i$.

    \item \textbf{Blocker source edges.}
    For each variable gadget $i$, create a dummy edge
    \[
        \eta_i \coloneqq (\bar d_i,\bar c_i),
    \]
    with nominal weight $1$. Add a zero-weight connector from $\bar c_i$ to the tail of the blocking edge $h_i$. Thus, each blocker source edge can route one unit of perturbation mass to the corresponding blocking edge. All edges not specified as source edges have nominal weight $0$.
\end{itemize}

We set the local diffusion budget to $\varepsilon\coloneqq 1$; by rescaling all edge weights, the same construction applies to any fixed positive budget. The dummy vertices $d_j,\bar d_i$ are unreachable from $s$, and the only incoming edges to $c_j,\bar c_i$ are the source edges $\alpha_j,\eta_i$, respectively. Hence the vertices $c_j,\bar c_i$ are also unreachable from $s$, and no directed $s$--$t$ path can use a source edge $\alpha_j,\eta_i$ or any connector edge leaving $c_j$ or $\bar c_i$. Therefore, every directed $s$--$t$ path in $G'$ remains on the variable backbone, chooses exactly one branch from each pair $(T_i,F_i)$, and hence encodes a truth assignment. Conversely, every truth assignment $A$ defines a canonical directed $s$--$t$ path, denoted by $P(A)$.

\begin{figure}[t]
    \centering
    \includegraphics[width=0.95\linewidth]{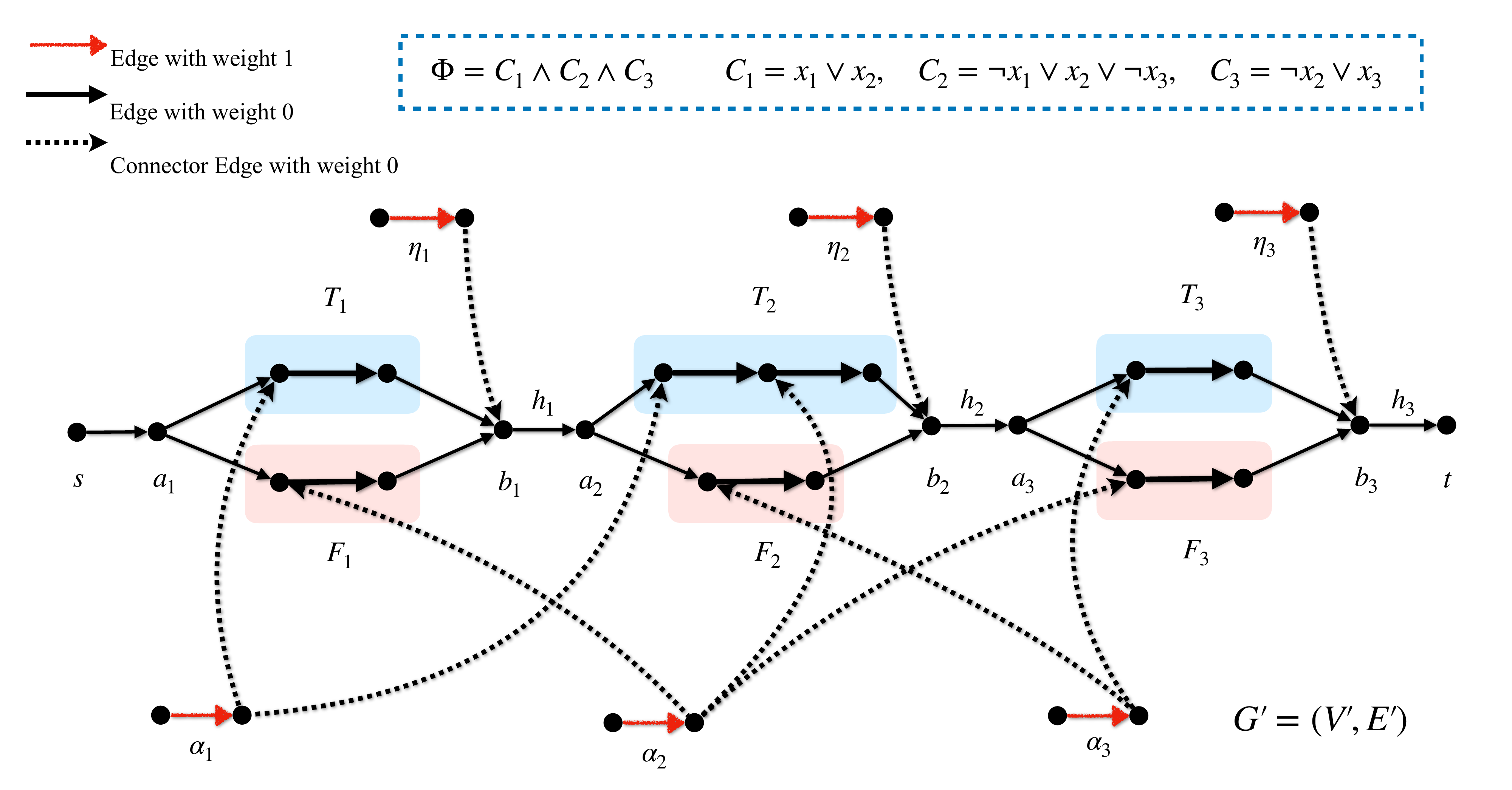}
    \caption{Example illustrating the reduction from MinSAT to Diff-RSP under $\mathcal D^{\mathrm{L},\infty}(\varepsilon)$. The formula is $\Phi=C_1\wedge C_2\wedge C_3$, with $C_1=x_1\vee x_2$, $C_2=\neg x_1\vee x_2\vee \neg x_3$, and $C_3=\neg x_2\vee x_3$. The graph $G'=(V',E')$ contains one variable gadget for each variable $x_i$, with true and false branches $T_i$ and $F_i$. Clause source edges $\alpha_j$ connect to the literal-occurrence edges of the corresponding clauses, and blocker source edges $\eta_i$ connect to the blocking edges $h_i$. Edges drawn in red have nominal weight $1$; all other edges have nominal weight $0$.}
    \label{fig:rsp-linf-example}
\end{figure}

\begin{example}[Illustration of the local-budget reduction]
\label{ex:rsp-linf-construction}
Figure~\ref{fig:rsp-linf-example} illustrates the construction on the MinSAT instance
\[
    \Phi = C_1\wedge C_2\wedge C_3,
    \qquad
    C_1=x_1\vee x_2,\quad
    C_2=\neg x_1\vee x_2\vee \neg x_3,\quad
    C_3=\neg x_2\vee x_3.
\]
The three variable gadgets are arranged in series, and choosing one branch from each pair $(T_i,F_i)$ encodes a truth assignment. The blocker source edges $\eta_i$ provide the fixed baseline contribution through the blocking edges $h_i$. Each clause source edge $\alpha_j$ connects to the occurrence edges corresponding to the literals in $C_j$; for instance, $\alpha_2$ connects to the branches representing $\neg x_1$, $x_2$, and $\neg x_3$. Thus, satisfied clauses create one additional unit contribution to the robust cost through a chosen branch, whereas unsatisfied clauses cannot create an additional contribution beyond the blocker baseline. \hfill $\clubsuit$
\end{example}

We now formalize this reduction. Let
\[
    \ell_\Phi \coloneqq \sum_{j=1}^m |C_j|
\]
denote the total number of literal occurrences in $\Phi$. The next proposition shows that, for every directed $s$--$t$ path $P$ in the constructed graph $G'$, the robust value under $\mathcal D^{\mathrm{L},\infty}(1)$ is exactly the fixed baseline $n$ plus $\mathrm{sat}(A(P))$, where $A(P)$ is the truth assignment encoded by $P$. Thus, minimizing the robust cost of a path in $G'$ is equivalent to choosing an assignment that minimizes the number of satisfied clauses, thereby establishing assertion~(iii) of Theorem~\ref{thm:diff-rsp}.

\begin{proposition}[Reduction from MinSAT to Diff-RSP under $\mathcal D^{\mathrm{L},\infty}(1)$]
\label{prop:rsp-linf-hardness}
Given an instance of MinSAT with formula
\[
    \Phi=C_1\wedge\cdots\wedge C_m
\]
over variables $x_1,\ldots,x_n$, one can construct in $O(n+\ell_\Phi)$ time a Diff-RSP instance on a directed graph $G'=(V',E')$ under $\mathcal D^{\mathrm{L},\infty}(1)$. Moreover, every directed $s$--$t$ path $P$ in $G'$ induces a truth assignment $A(P)$, and its robust value satisfies
\[
    \max_{(\Delta^+,\Delta^-)\in \mathcal D^{\mathrm{L},\infty}(1)}
    f^\top\bigl(w+\Delta^+-\Delta^-\bigr)
    =
    n+\mathrm{sat}(A(P)),
\]
where $f$ is the incidence vector of $P$. Conversely, for every truth assignment $A$, the corresponding canonical path $P(A)$ has robust value $n+\mathrm{sat}(A)$. In particular, MinSAT $\le_p$ Diff-RSP under $\mathcal D^{\mathrm{L},\infty}(\varepsilon)$, and Diff-RSP is \textsf{NP}-hard under $\mathcal D^{\mathrm{L},\infty}(\varepsilon)$.
\end{proposition}

Proposition~\ref{prop:rsp-linf-hardness} shows that the long-term local-budget regime can encode a clause-counting objective exactly. The blocker source edges produce the assignment-independent baseline $n$, while the clause source edges add precisely the assignment-dependent term $\mathrm{sat}(A)$. Thus, the absence of a global transportation charge does not make the problem easier; instead, it allows many local diffusion channels to be activated simultaneously, turning robust shortest path selection into the problem of choosing a truth assignment with the fewest satisfied clauses.


\subsection{\textsf{NP}-hardness under long-term global budget}
\label{subsec:rsp-np-global}

The long-term global $\ell_1$-budget regime becomes hard for the same basic reason: the adversary can move perturbation mass over multiple steps and concentrate it on the particular edges used by the chosen path. As in the local-budget regime, this breaks the short-term ``worst-case graph'' reduction because the adversarial effect is no longer purely one-step and node-local. The difference is that, under a global budget, every transported unit consumes part of a single shared $\ell_1$ budget. Thus, multi-step diffusion creates not only a reachability effect, but also a transportation-cost effect: mass originating several hops away can support the worst-case increase on a path edge, but only by paying for the route through which it travels. The key implication is that the adversary's impact on a candidate $s$--$t$ path is no longer determined by local neighborhoods independently at each node; it depends on how the path is situated in the directed topology, i.e., on which nodes lie in its neighborhood and can participate in multi-step propagation at sufficiently low transportation cost.

This global-budget dependence on the \emph{surroundings} of a path is naturally captured by the \emph{Most Secluded Path} problem, which asks for an $s$--$t$ path whose exposure to its neighborhood is minimal.

\begin{definition}[Most Secluded Path]
\label{def:most-secluded-path}
Let $G=( V, E)$ be a directed graph with terminals $s,t\in V$. For a node set $\mathcal S\subseteq  V$, define its closed \emph{out-neighborhood} by
\[
    N[\mathcal S]
    \;\coloneqq\;
    \mathcal S \cup \bigl\{ v\in  V : \exists\, u\in \mathcal S \text{ such that } (u,v)\in E \bigr\}.
\]
For an $s$--$t$ directed path $Q$, let $ V(Q)\subseteq V$ denote its node set and define its \emph{exposure} by
\[
    \mathrm{exp}(Q) \;\coloneqq\; \bigl|N[ V(Q)]\bigr|.
\]
The \emph{Most Secluded Path} problem asks to find an $s$--$t$ directed path $Q$ minimizing $\mathrm{exp}(Q)$.
\end{definition}

Most Secluded Path is computationally intractable in graph classes closely aligned with our reduction. In particular, Chechik, Johnson, Parter, and Peleg~\cite[Corollary~3.1]{chechik2017secluded} show that Most Secluded Path is \textsf{NP}-hard even on directed graphs of maximum degree four.\footnote{We cite this result only as a hardness source; the reduction below is explicit and tailored to the diffusion model.}

The connection to long-term global diffusion is now intuitive. If a chosen $s$--$t$ path $Q$ is highly exposed, then many nodes in its closed out-neighborhood $N[V(Q)]$ can supply perturbation mass that, under multi-step propagation, can be routed onto the corresponding path in the constructed instance and increase its robust cost; if $Q$ is secluded, then far fewer nodes can do so. The global $\ell_1$ budget turns this reachability into a transportation effect: exposed nodes can reach the canonical path cheaply, whereas nodes outside the exposure set can contribute only through routes whose transportation cost is prohibitively large. The reduction below makes this correspondence precise by constructing an instance in which, for every $s$--$t$ path $Q$ in the original graph, the worst-case long-term global-budget cost of the corresponding path lies between $|N[V(Q)]|$ and $|N[V(Q)]|+\tfrac{1}{2}$. The additive gap below $\tfrac{1}{2}$ is the key feature: since $|N[V(Q)]|$ is integer-valued, it is sufficient to recover the threshold structure of Most Secluded Path, thereby establishing assertion~(iv) of Theorem~\ref{thm:diff-rsp}. We now give the construction.

\medskip
\noindent\textbf{Construction of $G''$.}
The construction replaces each original node $v\in V$ by a directed \emph{edge-chain} gadget $E_v$ together with a private entry node $c_v$, a single \emph{dummy edge} of nominal weight $1$, and a zero-weight entry edge into the gadget. It also introduces two types of inter-gadget edges. \emph{Anchor} edges preserve the original path structure of $G$, while \emph{connector} edges encode the closed out-neighborhood relation used in the definition of path exposure in Most Secluded Path (Definition~\ref{def:most-secluded-path}). The chain lengths and connector geometry are chosen to encode transportation cost under the global budget: reaching the canonical path through a connector edge is cheap, whereas reaching it through the anchor structure requires traversing a long chain and is therefore expensive. Figure~\ref{fig:rsp-construction-Gdoubleprime} illustrates the edge-chain gadget, the anchor edges, and the connector edges. The graph $G''$ is constructed as follows:

\begin{figure}[t]
    \centering
    \includegraphics[width=0.95\linewidth]{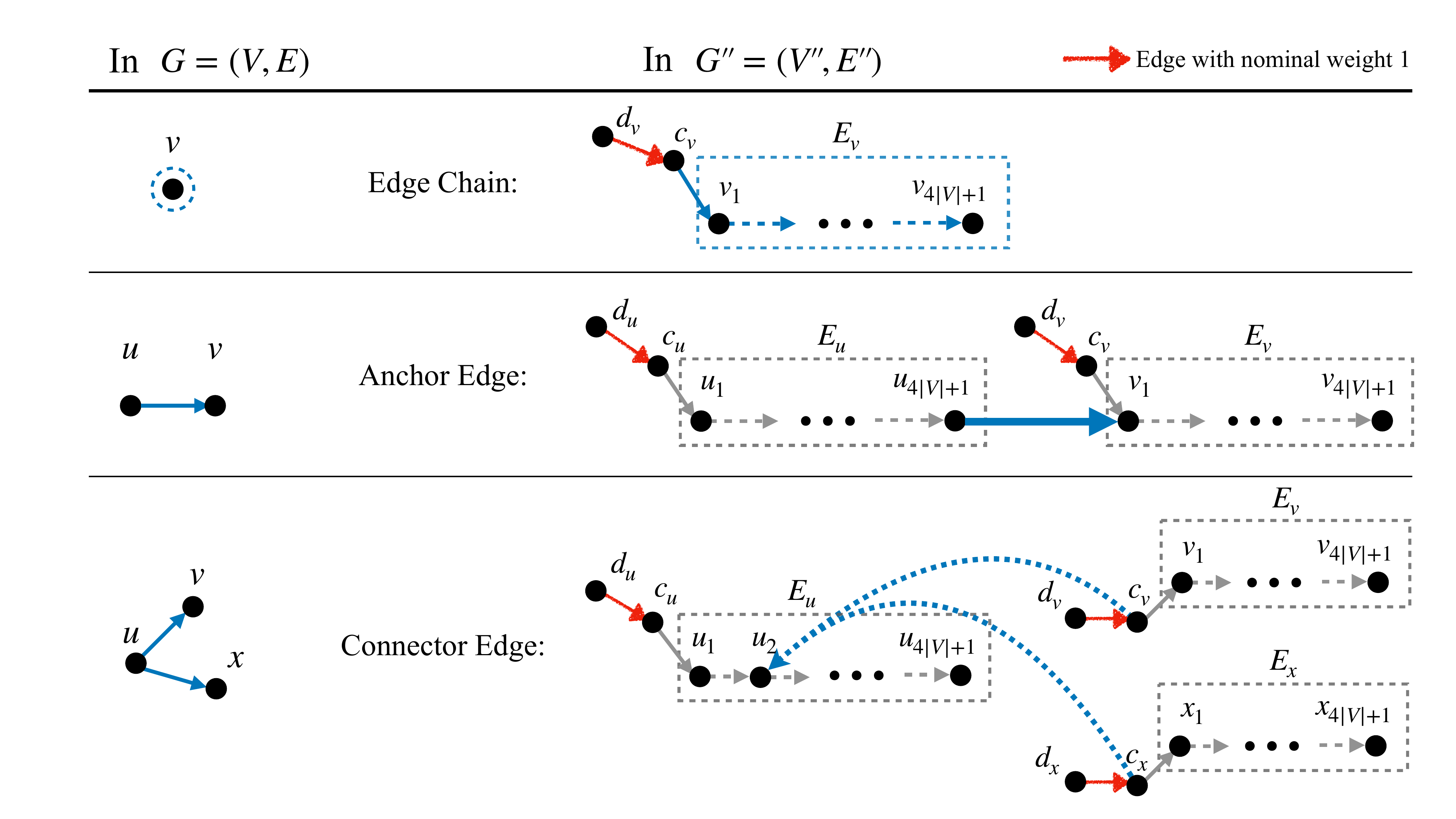}
    \caption{Construction of $G''=(V'',E'')$ from $G=(V,E)$. Top: each node $v\in V$ is replaced by a chain gadget $E_v$ of length $4|V|$, together with a dummy edge $\alpha_v=(d_v,c_v)$ of nominal weight $1$ and an entry edge $\beta_v=(c_v,v_1)$. Middle: each original arc $(u,v)\in E$ induces an anchor edge $a_{u,v}$ from the end of $E_u$ to the beginning of $E_v$ (solid blue). Bottom: each relation $v\in N[u]\setminus\{u\}$ induces a connector edge $c_{v,u}=(c_v,u_2)$ (dashed blue), allowing the unit from the gadget of $v$ to enter the gadget $E_u$ near its beginning. All edges other than the dummy edges have nominal weight $0$.}
    \label{fig:rsp-construction-Gdoubleprime}
\end{figure}

\begin{itemize}
    \item \textbf{Edge-chain gadget.}
    For each $v\in V$, create a directed chain
    \[
        E_v \coloneqq \{e_v^i:\, i=1,\ldots,4|V|\}\subseteq E''.
    \]
    The length of this chain is $4|V|$. Denote the vertices in this chain by $v_1,\ldots,v_{4|V|+1}$, so that $e_v^i=(v_i,v_{i+1})$, for $i=1,\ldots,4|V|$.

    \item \textbf{Dummy and entry edges.}
    For each $v\in V$ we add a dummy vertex $d_v$, an entry vertex $c_v$, a dummy edge
    \[
        \alpha_v \coloneqq (d_v,c_v),
    \]
    and an entry edge
    \[
        \beta_v \coloneqq (c_v,v_1).
    \]
    Each dummy edge $\alpha_v$ has nominal weight $1$, and every other edge in $G''$, including each entry edge $\beta_v$, has nominal weight $0$.

    \item \textbf{Anchor edges.}
    For each $(u,v)\in E$, create an anchor edge
    \[
        a_{u,v} \coloneqq (u_{4|V|+1},\,v_1).
    \]
    Thus, a canonical path moves from the end of $E_u$ to the beginning of $E_v$ exactly when $(u,v)\in E$.

    \item \textbf{Connector edges.}
    For each $u\in V$ and each $v\in N[u]\setminus\{u\}$, create a connector edge
    \[
        c_{v,u} \coloneqq (c_v,\,u_2).
    \]
    Equivalently, for every arc $(u,v)\in E$, the unit at the entry vertex $c_v$ can enter the gadget $E_u$ near its beginning, at the vertex $u_2$.
\end{itemize}

We set the diffusion budget to $\varepsilon\coloneqq 4|V|$. Since the dummy edges are the only edges with positive nominal weight, each gadget can emit at most one unit of perturbation mass. The construction is arranged so that every node in the exposure set $N[V(Q)]$ can route its unit onto the canonical path at transportation cost $4$. Thus realizing all exposed contributions uses total budget $4|N[V(Q)]|$ and leaves remaining budget $4\bigl(|V|-|N[V(Q)]|\bigr)$. By contrast, any contribution from outside the exposure set has no connector directly into an on-path gadget. Hence, before it can first reach the canonical path, it must use a non-dummy transfer edge out of its entry node and traverse at least $4|V|-1$ chain edges outside the path, requiring transportation cost at least $8|V|+2$ per unit. Consequently, the remaining budget can increase the robust value by less than
\[
    \frac{4(|V|-|N[V(Q)]|)}{8|V|+2}<\frac{1}{2}.
\]
It follows that, for every path $Q$ in the original graph, the robust value of the corresponding canonical path in $G''$ lies between $|N[V(Q)]|$ and $|N[V(Q)]|+\tfrac{1}{2}$.

\begin{figure}[t]
    \centering
    \includegraphics[width=0.9\linewidth]{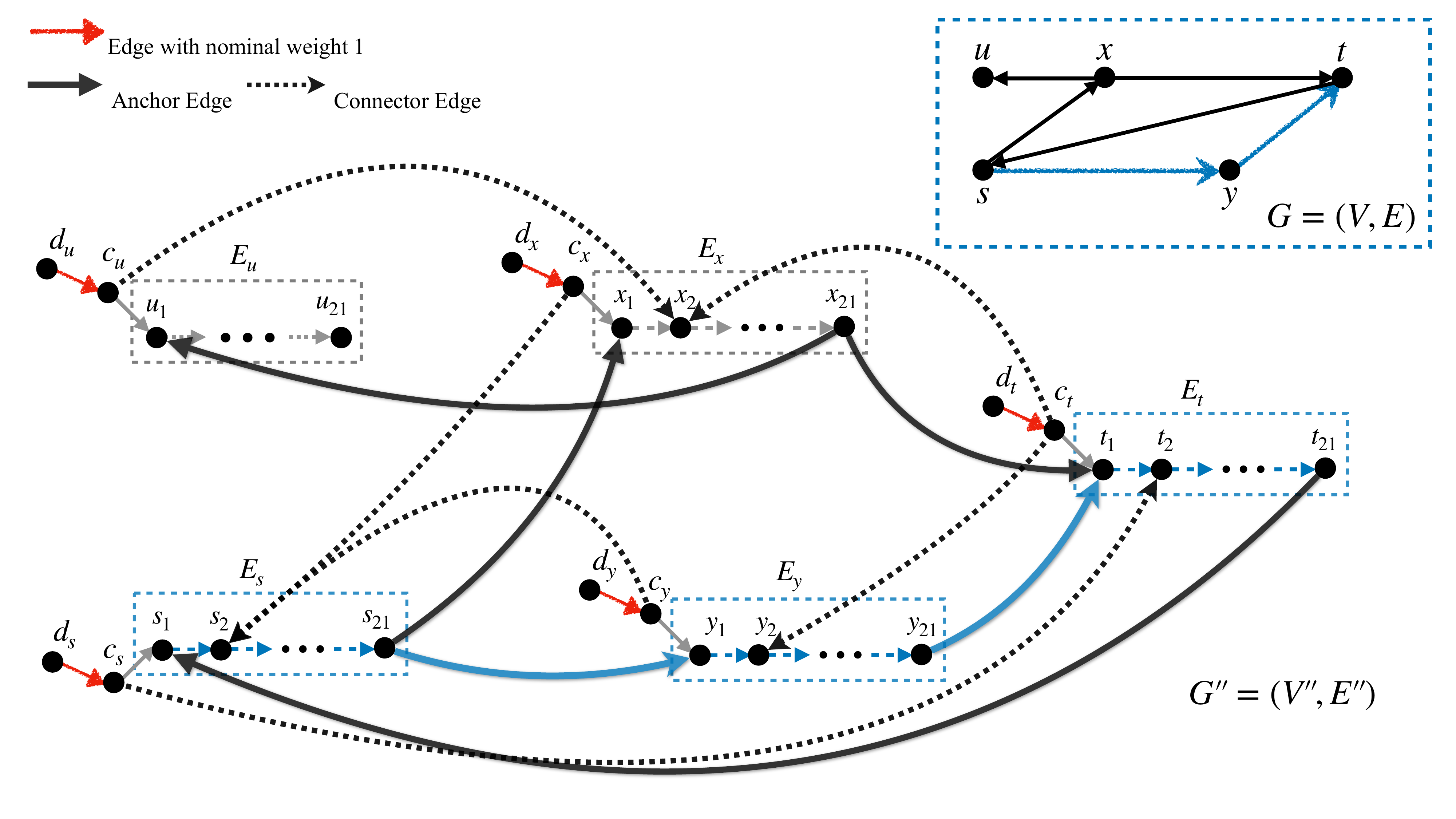}
    \caption{Example illustrating the reduction from Most Secluded Path to Diff-RSP under $\mathcal D^{\mathrm{L},1}(\varepsilon)$. Top right: a directed graph $G=(V,E)$ with highlighted path $Q:s\to y\to t$ (blue). Bottom: the corresponding constructed graph $G''=(V'',E'')$. Solid edges are anchor edges and dotted edges are connector edges. The highlighted blue path is the canonical path $P(Q)$ corresponding to $Q$.}
    \label{fig:rsp-l1-example}
\end{figure}

\begin{example}[Illustration of the global-budget reduction]
\label{ex:rsp-l1-construction}
Figure~\ref{fig:rsp-l1-example} illustrates the construction on a small directed graph $G=(V,E)$. The highlighted path in $G$ is $Q:s\to y\to t$, and the highlighted blue path in the constructed graph $G''=(V'',E'')$ is the corresponding canonical path $P(Q)$. The example illustrates the key feature of the construction: connector routes allow nodes in the exposure set $N[V(Q)]$ to reach the canonical path cheaply, whereas nodes outside the exposure set can reach it only through substantially longer anchor-based routes. \hfill $\clubsuit$
\end{example}

We now formalize this reduction. The next proposition associates each directed $s$--$t$ path $Q$ in $G$ with a canonical directed $s'$--$t'$ path $P(Q)$ in the constructed graph $G''$, and proves lower and upper bounds on the robust value of $P(Q)$ under $\mathcal D^{\mathrm{L},1}(4|V|)$ that differ by less than $\tfrac{1}{2}$. Together, these bounds show that the global-budget construction tracks the exposure of $Q$ closely enough to recover the threshold structure of Most Secluded Path, thereby establishing assertion~(iv) of Theorem~\ref{thm:diff-rsp}.

\begin{proposition}[Reduction from Most Secluded Path to Diff-RSP under $\mathcal D^{\mathrm{L},1}(4|V|)$]
\label{prop:rsp-l1-hardness}
Given an instance of Most Secluded Path on a directed graph $G=(V,E)$ with terminals $s,t\in V$, one can construct in $O(|V|^2+|E|)$ time a Diff-RSP instance on a directed graph $G''=(V'',E'')$ under $\mathcal D^{\mathrm{L},1}(4|V|)$. Moreover, for every directed $s$--$t$ path $Q$ in $G$, the corresponding canonical directed $s'$--$t'$ path $P(Q)$ in $G''$ has robust value satisfying
\[
    |N[V(Q)]|
    \;\le\;
    \max_{(\Delta^+,\Delta^-)\in \mathcal D^{\mathrm{L},1}(4|V|)}
    f^\top\bigl(w+\Delta^+-\Delta^-\bigr)
    \;<\;
    |N[V(Q)]|+\tfrac{1}{2},
\]
where $f$ is the incidence vector of $P(Q)$. In particular, Most Secluded Path $\le_p$ Diff-RSP under $\mathcal D^{\mathrm{L},1}(\varepsilon)$, and Diff-RSP is \textsf{NP}-hard under $\mathcal D^{\mathrm{L},1}(\varepsilon)$.
\end{proposition}

Proposition~\ref{prop:rsp-l1-hardness} shows that the global-budget construction encodes the minimum-exposure objective of Most Secluded Path up to an additive gap of less than $\tfrac{1}{2}$. Because $|N[V(Q)]|$ is always an integer, this is sufficient for a threshold reduction. The longer gadgets are what make the global $\ell_1$ budget act as a transportation budget: their length is scaled to $4|V|$ so that connector routes into the canonical path remain cheap, whereas anchor-based routes from outside the exposure set become expensive.


\section{Diffusion-Robust Traveling Salesman Problem}
\label{sec:tsp}

The traveling salesman problem (TSP) provides a useful contrast with the shortest-path problem studied in Section~\ref{sec:shortest-path}. In the nominal setting, TSP is already \textsf{NP}-hard, and therefore the relevant question is not whether diffusion creates hardness from an otherwise tractable problem, but how the structure of diffusion changes the robust counterpart relative to ordinary TSP. We show that three of the four uncertainty sets preserve the ordinary TSP structure in a strong sense: under $\mathcal D^{\mathrm{S},\infty}(\varepsilon)$, $\mathcal D^{\mathrm{S},1}(\varepsilon)$, and $\mathcal D^{\mathrm{L},1}(\varepsilon)$, Diff-RTSP is polynomial-time equivalent to ordinary TSP. The underlying reason is that, once a Hamiltonian cycle is fixed, every vertex has exactly one incoming tour edge and one outgoing tour edge. This structure makes the fixed-tour worst-case value collapse either to the cost of the same tour under modified edge weights or to a capped expression depending only on its nominal cost. The remaining regime, $\mathcal D^{\mathrm{L},\infty}(\varepsilon)$, is qualitatively different: multi-step propagation combined with local edgewise caps creates a capacitated diffusion problem inside the evaluation of a fixed tour, so this case does not collapse by the same fixed-tour argument to ordinary TSP. Nevertheless, we derive a computable upper bound on its optimal robust value using two ordinary TSP instances.

Let $G=(V,E)$ be a \emph{complete} directed graph with nonnegative edge weights $w\in\mathbb R_+^E$. A Hamiltonian cycle is a directed cycle that visits every vertex in $V$ exactly once and returns to its starting vertex. We denote by $\mathcal H$ the set of Hamiltonian cycles in $G$. For $H\in\mathcal H$, write
\[
    W(H):=\sum_{e\in H}w_e,
\]
and let
\[
    \mathrm{OPT}(w)
    :=
    \min_{H\in\mathcal H}W(H)
\]
denote the optimal value of ordinary TSP with weights $w$. Given a diffusion uncertainty set $\mathcal D(\varepsilon)$, the diffusion-robust traveling salesman problem, abbreviated Diff-RTSP, is
\begin{equation}
\label{eq:diff-rtsp}
    \mathrm{OPT}_{\mathrm{RTSP}}(w,\mathcal D(\varepsilon))
    :=
    \min_{H\in\mathcal H}
    \max_{(\Delta^+,\Delta^-)\in\mathcal D(\varepsilon)}
    \sum_{e\in H}\bigl(w_e+\Delta_e^+-\Delta_e^-\bigr).
\end{equation}
Equivalently, if $f_H\in\{0,1\}^E$ denotes the incidence vector of the tour $H$, then the objective in~\eqref{eq:diff-rtsp} can be written as $f_H^\top(w+\Delta^+-\Delta^-)$.

We now state the main reduction result for the three regimes that reduce to ordinary TSP. As in the shortest-path case, the robust problem is at least as hard as ordinary TSP, because ordinary TSP is recovered as the special case $\varepsilon=0$. The substantive part is the converse direction: under each of $\mathcal D^{\mathrm{S},\infty}(\varepsilon)$, $\mathcal D^{\mathrm{S},1}(\varepsilon)$, and $\mathcal D^{\mathrm{L},1}(\varepsilon)$, Diff-RTSP admits a polynomial-time reduction to ordinary TSP.

\begin{theorem}[Reduction of Diff-RTSP to ordinary TSP]
\label{thm:diff-rtsp}
For each uncertainty set
\[
    \mathcal D(\varepsilon)
    \in
    \left\{
    \mathcal D^{\mathrm{S},\infty}(\varepsilon),
    \mathcal D^{\mathrm{S},1}(\varepsilon),
    \mathcal D^{\mathrm{L},1}(\varepsilon)
    \right\},
\]
Diff-RTSP is polynomial-time equivalent to ordinary TSP. In particular, ordinary TSP is recovered as the special case $\varepsilon=0$, and for each of the three uncertainty sets above, Diff-RTSP can be reduced in polynomial time to ordinary TSP.
\end{theorem}

The proof is based on fixed-tour characterizations. We fix an arbitrary Hamiltonian cycle $H$ and compute the value of the inner maximization in~\eqref{eq:diff-rtsp}. This mirrors the role of pathwise characterizations in the shortest-path section, but the organization is different. For Diff-RSP, the main distinction was between tractable and intractable regimes. For Diff-RTSP, the three regimes covered by Theorem~\ref{thm:diff-rtsp} are all reducible to ordinary TSP, and the relevant distinction is instead the form of the fixed-tour worst-case value. Under the short-term local budget $\mathcal D^{\mathrm{S},\infty}(\varepsilon)$, the worst-case value induces modified edge weights. Under the two global-budget regimes, $\mathcal D^{\mathrm{S},1}(\varepsilon)$ and $\mathcal D^{\mathrm{L},1}(\varepsilon)$, the worst-case value has a capped form that depends on the tour only through its nominal weight. We treat these two mechanisms first, and then discuss the remaining long-term local-budget regime separately.


\subsection{Modified edge weights under short-term local budget}
\label{subsec:tsp-short-infty}

We begin with the short-term local budget set $\mathcal D^{\mathrm{S},\infty}(\varepsilon)$. Among the four TSP regimes, this is the only one in which the worst-case value of a fixed tour depends on the local edge structure of the tour beyond its nominal weight. The reason is that, under the short-term propagation constraint, mass removed from incoming edges into a vertex $u$ can only be redistributed to edges leaving $u$. Once a Hamiltonian cycle is fixed, there is a unique outgoing tour edge from $u$, so any mass that the adversary wants to add to the tour at $u$ must be placed on this edge. The useful supply at $u$, however, comes only from incoming edges that are not used by the tour: removing mass from the unique incoming tour edge would directly decrease the tour cost. Thus, the adversary's gain at each vertex is determined by the off-tour incoming supply at that vertex, capped by the local budget $\varepsilon$.

This observation leads to the same type of preprocessing used for Diff-RSP under the short-term local budget. For each vertex $u\in V$, define the one-step transfer bound
\[
    T_u \coloneqq \sum_{a\in E_{\mathrm{in}}(u)} \min\{\varepsilon,w_a\}.
\]
This is the maximum perturbation mass that can be collected from the incoming neighborhood of $u$ under short-term diffusion. For each edge $e=(i,u)\in E$, define the surcharge
\[
    \chi_e
    \coloneqq
    \min\left\{
        \varepsilon,\,
        T_u-\min\{\varepsilon,w_e\}
    \right\} \;=\; \min\left\{\varepsilon,\; \sum_{e^\prime \in  E_{\mathrm{in}}(u)\setminus\{e\}} \min\{\varepsilon,w_{e^\prime}\}\right\}.
\]
Finally, let
\[
    w_e^{\mathrm{wc}} \coloneqq w_e+\chi_e,
    \qquad e\in E,
\]
and write $w^{\mathrm{wc}}\in\mathbb R_+^E$ for the resulting edge-weight vector of the precomputed worst-case graph. The quantity $\chi_e$ is the maximum gain that can be generated at the head vertex $u$ when $e$ is the tour edge entering $u$: all other incoming edges into $u$ are off-tour and can supply at most $T_u-\min\{\varepsilon,w_e\}$ units of mass in total, while the unique outgoing tour edge from $u$ can receive at most $\varepsilon$. We charge this gain to the incoming edge $e$, which allows the fixed-tour worst-case value to be written using the precomputed weights $w^{\mathrm{wc}}$.

\begin{proposition}[Reduction to ordinary TSP under short-term local budget]
\label{prop:tsp-short-infty}
The weights $w^{\mathrm{wc}}$ can be computed in $O(|E|+|V|)$ time. Moreover, for every Hamiltonian cycle $H\in\mathcal H$,
\[
    \max_{(\Delta^+,\Delta^-)\in\mathcal D^{\mathrm{S},\infty}(\varepsilon)}
    \sum_{e\in H}\bigl(w_e+\Delta_e^+-\Delta_e^-\bigr)
    =
    \sum_{e\in H}w_e^{\mathrm{wc}}.
\]
Thus, under $\mathcal D^{\mathrm{S},\infty}(\varepsilon)$, Diff-RTSP reduces to one ordinary TSP computation with edge weights $w^{\mathrm{wc}}$.
\end{proposition}

Proposition~\ref{prop:tsp-short-infty} shows that the short-term local budget preserves the TSP structure after a deterministic edge-weight transformation. We next turn to the two global-budget regimes, where no such local edge modification is needed: the worst-case value of a fixed tour takes a capped form depending only on its nominal weight.


\subsection{Capped fixed-tour values under global budgets}
\label{subsec:tsp-global}

We now turn to the two global-budget regimes, $\mathcal D^{\mathrm{S},1}(\varepsilon)$ and $\mathcal D^{\mathrm{L},1}(\varepsilon)$. In these regimes, the robust value of a fixed tour no longer requires the local edgewise surcharges used in the short-term local-budget case. Instead, conservation and the $\ell_1$ budget impose a scalar upper bound $\varepsilon/2$ on the net additional mass that can be added to any tour. A second limitation is the total mass outside the tour: the adversary cannot move more nominal mass into the tour than is initially available on edges not used by the tour. These two bounds lead to a capped fixed-tour value depending only on the nominal weight of the tour.

For a Hamiltonian cycle $H\in\mathcal H$, recall that $W(H)=\sum_{e\in H}w_e$, and define
\[
    S:=\sum_{e\in E}w_e.
\]
Thus $S-W(H)$ is the total nominal mass on edges outside the tour. In this notation, the two relevant upper bounds are transparent: the $\ell_1$ budget limits the net gain to $\varepsilon/2$, while conservation of total post-diffusion mass limits the total tour cost to $S$. The next proposition shows that these bounds are tight, both in the short-term and long-term global-budget regimes.

\begin{proposition}[Global-budget regimes]
\label{prop:tsp-global-budget}
Let $\mathcal D(\varepsilon)$ be either $\mathcal D^{\mathrm{S},1}(\varepsilon)$ or $\mathcal D^{\mathrm{L},1}(\varepsilon)$. Then, for every Hamiltonian cycle $H\in\mathcal H$,
\[
    \max_{(\Delta^+,\Delta^-)\in\mathcal D(\varepsilon)}
    \sum_{e\in H}\bigl(w_e+\Delta_e^+-\Delta_e^-\bigr)
    =
    \min\left\{
        W(H)+\frac{\varepsilon}{2},
        S
    \right\}.
\]
Consequently, under either global-budget uncertainty set, Diff-RTSP reduces to ordinary TSP with the original weights $w$. More precisely,
\[
    \mathrm{OPT}_{\mathrm{RTSP}}(w,\mathcal D(\varepsilon))
    =
    \min\left\{
        \mathrm{OPT}(w)+\frac{\varepsilon}{2},
        S
    \right\}.
\]
\end{proposition}

Proposition~\ref{prop:tsp-global-budget} completes the proof of Theorem~\ref{thm:diff-rtsp} for the two global-budget regimes. Together with Proposition~\ref{prop:tsp-short-infty}, it shows that three of the four TSP regimes reduce to ordinary TSP through explicit fixed-tour formulas. The remaining regime, $\mathcal D^{\mathrm{L},\infty}(\varepsilon)$, is different because multi-step propagation and local edgewise caps interact: evaluating a fixed tour involves a capacitated diffusion problem, in the sense that perturbation mass may be redistributed across vertices before being added to the tour, while every edge used to receive or relay this mass remains subject to the same componentwise local budget. We discuss this case next.


\subsection{The long-term local-budget regime}
\label{subsec:tsp-long-infty}

The long-term local-budget set $\mathcal D^{\mathrm{L},\infty}(\varepsilon)$ combines multi-step propagation with componentwise local caps. This combination changes the fixed-tour evaluation problem. Under long-term propagation, perturbation mass can be redistributed across vertices before being added to the tour. However, this redistribution is itself capacity-constrained: every edge used to receive or relay mass can receive at most $\varepsilon$ units and forward at most $\varepsilon$ units. Thus, even in a complete graph, the aggregate amount of capped off-tour mass need not be simultaneously deliverable to the tour edges. In this subsection, we derive a capped upper bound that can be computed from two ordinary TSP instances.

For comparison with the preceding regimes, define
\[
    c_e:=\min\{\varepsilon,w_e\},
    \qquad \text{for }e\in E,
\]
and let
\[
    C:=\sum_{e\in E}c_e.
\]
The quantity $c_e$ is the maximum net amount of nominal mass that can be drained from edge $e$ under the local budget. Therefore, for a fixed Hamiltonian cycle $H$, the total locally capped nominal mass outside the tour is
\[
    \sum_{e\notin H}c_e
    =
    C-\sum_{e\in H}c_e.
\]
At the same time, each tour edge can receive at most $\varepsilon$ units of perturbation mass, so the total gain on the $n$ tour edges is also bounded by $n\varepsilon$. These two observations yield a simple computable upper bound, obtained from two ordinary TSP computations.

\begin{lemma}[Capped upper bound under long-term local budget]
\label{lem:tsp-long-infty-upper}
For every Hamiltonian cycle $H\in\mathcal H$,
\[
    \max_{(\Delta^+,\Delta^-)\in\mathcal D^{\mathrm{L},\infty}(\varepsilon)}
    \sum_{e\in H}\bigl(w_e+\Delta_e^+-\Delta_e^-\bigr)
    \le
    \min\left\{
        W(H)+n\varepsilon,\,
        C+\sum_{e\in H}(w_e-c_e)
    \right\}.
\]
Consequently,
\[
    \mathrm{OPT}_{\mathrm{RTSP}}\bigl(w,\mathcal D^{\mathrm{L},\infty}(\varepsilon)\bigr)
    \le
    \min\left\{
        \mathrm{OPT}(w)+n\varepsilon,\,
        C+\mathrm{OPT}(w-c)
    \right\}.
\]
The quantities $c$, $w-c$, and $C$ can be computed in $O(|E|)$ time.
\end{lemma}

The upper bound in Lemma~\ref{lem:tsp-long-infty-upper} deliberately uses only the two aggregate bounds that lead to ordinary TSP computations: the total receiving capacity of the tour and the total locally capped mass outside the tour. Incorporating routing bottlenecks can lead to sharper bounds, but such bounds generally depend on the interaction between the tour and the capacity-constrained diffusion problem and no longer have the same immediate reduction to ordinary TSP. In particular, the bound in Lemma~\ref{lem:tsp-long-infty-upper} is not, in general, an exact fixed-tour formula. The following example shows that the capped upper bound can be strict, even on a complete directed graph.

\begin{example}[Strictness of the capped upper bound]
\label{ex:tsp-long-infty-strict}
Consider the complete directed graph on $V=\{0,1,2,3\}$, and fix the Hamiltonian cycle
\[
    H=\{(0,1),(1,2),(2,3),(3,0)\}.
\]
\begin{figure}[t]
    \centering
    \includegraphics[width=0.9\linewidth]{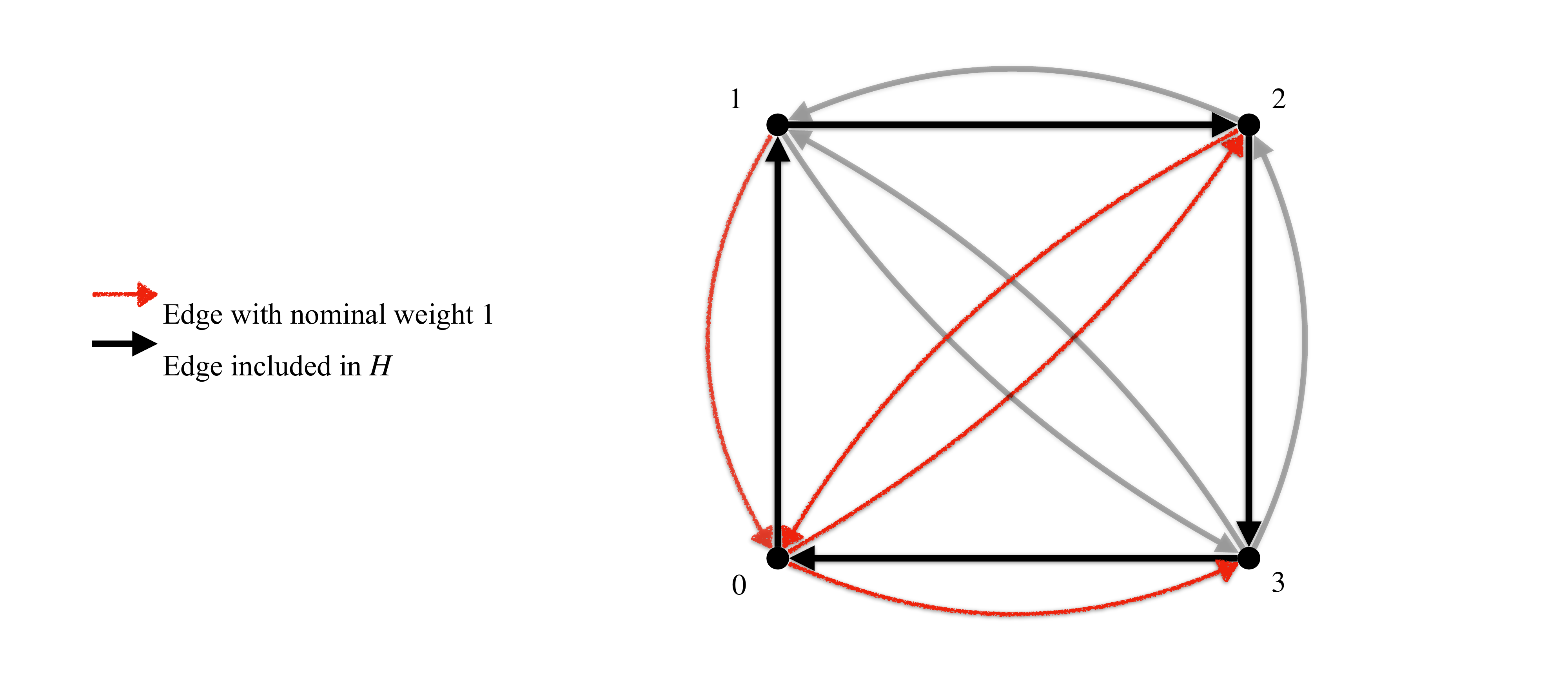}
    \caption{Visualization of the graph in Example \ref{ex:tsp-long-infty-strict}}
    \label{fig:tsp-strict-capped-upper-bound}
\end{figure}

Let $\varepsilon=1$. Set the weights of all tour edges to zero, set
\[
    w_{02}=w_{03}=w_{10}=w_{20}=1,
\]
and set all remaining off-tour edge weights to zero, as shown in Figure~\ref{fig:tsp-strict-capped-upper-bound}. Then $W(H)=0$ and
\[
    \sum_{e\notin H}\min\{1,w_e\}=4,
    \qquad
    n\varepsilon=4.
\]
Thus Lemma~\ref{lem:tsp-long-infty-upper} gives the upper bound of four on the worst-case gain of $H$.

We show that this upper bound cannot be attained. Suppose, for contradiction, that a feasible diffusion attained a gain of four. Then each of the four tour edges would have to receive its full local budget and lose no mass, so $\Delta_e^+=1$ and $\Delta_e^-=0$ for every $e\in H$. Now consider conservation at vertices $1$, $2$, and $3$. Since the off-tour edges $(1,3)$, $(2,1)$, $(3,1)$, and $(3,2)$ have zero nominal weight, the long-term constraint implies $\Delta^-_{ij}\le \Delta^+_{ij}$ on each of these four edges. Hence conservation at vertex $1$ gives
\[
    \Delta_{12}^++\Delta^+_{10}+\Delta^+_{13}
    =
    \Delta^-_{21}+\Delta^-_{31}
    \le
    \Delta^+_{21}+\Delta^+_{31}.
\]
Using conservation at vertices $2$ and $3$, together with the local-budget bounds $\Delta^-_{02}\le 1$ and $\Delta^-_{03}\le 1$ and the full-gain conditions $\Delta^+_{23}=\Delta^+_{30}=1$, we obtain $\Delta^+_{20}+\Delta^+_{21}\le \Delta^+_{32}$ and $\Delta^+_{31}+\Delta^+_{32}\le \Delta^+_{13}$. Combining these inequalities,
\[
    \Delta^+_{21}+\Delta^+_{31}
    \le
    \Delta^+_{32}+\Delta^+_{31}
    \le
    \Delta^+_{13}.
\]
Since $\Delta_{12}^+=1$ because $(1,2)\in H$, the preceding displays imply
\[
    1+\Delta^+_{10}+\Delta^+_{13}\le \Delta^+_{13},
\]
a contradiction. Thus the capped upper bound is strict for this fixed tour.  \hfill $\clubsuit$
\end{example}

Example~\ref{ex:tsp-long-infty-strict} shows that, even in a complete directed graph, locally available off-tour mass need not be simultaneously deliverable to the tour edges. This capacity-constrained redistribution is the feature that distinguishes the long-term local-budget regime from the three regimes covered by Theorem~\ref{thm:diff-rtsp}. Nevertheless, the value in this regime can still be bracketed by ordinary-TSP quantities.

\begin{remark}[Ordinary-TSP bounds for the long-term local-budget regime]
\label{rem:tsp-long-infty-sandwich}
Since $\mathcal D^{\mathrm{S},\infty}(\varepsilon)\subseteq\mathcal D^{\mathrm{L},\infty}(\varepsilon)$, because $\Delta^-\le w$ implies $\Delta^-\le w+\Delta^+$, Proposition~\ref{prop:tsp-short-infty} gives the lower bound
\[
    \mathrm{OPT}(w^{\mathrm{wc}})
    \le
    \mathrm{OPT}_{\mathrm{RTSP}}\bigl(w,\mathcal D^{\mathrm{L},\infty}(\varepsilon)\bigr).
\]
Combining this with Lemma~\ref{lem:tsp-long-infty-upper}, we obtain
\[
    \mathrm{OPT}(w^{\mathrm{wc}})
    \le
    \mathrm{OPT}_{\mathrm{RTSP}}\bigl(w,\mathcal D^{\mathrm{L},\infty}(\varepsilon)\bigr)
    \le
    \min\left\{
        \mathrm{OPT}(w)+n\varepsilon,\,
        C+\mathrm{OPT}(w-c)
    \right\}.
\]
Thus, although the long-term local-budget regime does not generally collapse to an ordinary TSP objective, its optimal value is bracketed by quantities computable from ordinary TSP instances. \hfill $\clubsuit$
\end{remark}

\section*{Acknowledgements}

Liviu Aolaritei acknowledges support from the Swiss National Science Foundation through the Postdoc.Mobility Fellowship (grant agreement P500PT\_222215). Funded in part by the European Union (ERC-2022-SYG-OCEAN-101071601). Views and opinions expressed are however those of the author(s) only and do not necessarily reflect those of the European Union or the European Research Council Executive Agency. Neither the European Union nor the granting authority can be held responsible for them.
Paul Grigas acknowledges the support of the NSF AI Institute for Advances in Optimization Award 2112533.


\bibliographystyle{abbrvnat} 
\bibliography{bibfile.bib}

\appendix

\section{Proofs}

\subsection{Proofs for Section~\ref{sec:shortest-path}}

\subsubsection{Proof of Theorem~\ref{thm:diff-rsp}}

Assertions~(i) and (ii) follow from Proposition~\ref{prop:diff-rsp-short-term-algorithm}. Assertion~(iii) follows from Proposition~\ref{prop:rsp-linf-hardness}. Assertion~(iv) follows from Proposition~\ref{prop:rsp-l1-hardness}.

\subsubsection{Proof of Proposition~\ref{prop:diff-rsp-short-term-algorithm}}

\textbf{(i) Case $\mathcal D^{\mathrm{S},\infty}(\varepsilon)$.} We first establish correctness, and then derive the running-time bound.

\medskip
\noindent\emph{Step 1: Pathwise equivalence.}
Fix any $f\in\mathcal P_{s,t}$, and let the corresponding directed $s$--$t$ path be
\[
    s=v_0 \xrightarrow{e_1} v_1 \xrightarrow{e_2} \cdots \xrightarrow{e_k} v_k=t.
\]
We claim that
\begin{equation}
\label{eq:pathwise-identity-infty}
    \max_{(\Delta^+,\Delta^-)\in \mathcal D^{\mathrm{S},\infty}(\varepsilon)}
    f^\top\bigl(w+\Delta^+-\Delta^-\bigr)
    =
    (w^{\mathrm{wc}})^\top f + c_s.
\end{equation}

\smallskip
\noindent\emph{Upper bound.}
Fix any feasible $(\Delta^+,\Delta^-)\in\mathcal D^{\mathrm{S},\infty}(\varepsilon)$. Regrouping the path objective gives
\begin{equation}
\label{eq:rsp-path-regroup}
    f^\top\bigl(w+\Delta^+-\Delta^-\bigr)
    =
    \sum_{i=1}^k w_{e_i}
    + \Delta^+_{e_1}
    + \sum_{i=1}^{k-1}\bigl(\Delta^+_{e_{i+1}}-\Delta^-_{e_i}\bigr)
    - \Delta^-_{e_k}.
\end{equation}
At the source,
\[
    \Delta^+_{e_1}
    \le
    \sum_{e\in E_{\mathrm{in}}(s)} \Delta^-_e
    \le
    \sum_{e\in E_{\mathrm{in}}(s)} \min\{\varepsilon,w_e\}
    =
    T_s,
\]
and also $\Delta^+_{e_1}\le \varepsilon$, hence
\[
    \Delta^+_{e_1}\le \min\{\varepsilon,T_s\}=c_s.
\]
Now fix $i\in\{1,\dots,k-1\}$. By conservation at $v_i$,
\[
    \sum_{e\in E_{\mathrm{out}}(v_i)} \Delta^+_e
    =
    \sum_{e\in E_{\mathrm{in}}(v_i)} \Delta^-_e.
\]
Since $e_i$ and $e_{i+1}$ are the unique path edges entering and leaving $v_i$,
\[
\begin{aligned}
    \Delta^+_{e_{i+1}}-\Delta^-_{e_i}
    \le
    \sum_{e\in E_{\mathrm{out}}(v_i)} \Delta^+_e - \Delta^-_{e_i} 
    =
    \sum_{e\in E_{\mathrm{in}}(v_i)\setminus\{e_i\}} \Delta^-_e 
    \le
    \sum_{e\in E_{\mathrm{in}}(v_i)\setminus\{e_i\}} \min\{\varepsilon,w_e\}
    =
    T_{v_i}-\min\{\varepsilon,w_{e_i}\}.
\end{aligned}
\]
Moreover,
\[
    \Delta^+_{e_{i+1}}-\Delta^-_{e_i}\le \Delta^+_{e_{i+1}}\le \varepsilon.
\]
Therefore,
\begin{equation}
\label{eq:rsp-node-bound}
    \Delta^+_{e_{i+1}}-\Delta^-_{e_i}
    \le
    \min\bigl\{\varepsilon,\;T_{v_i}-\min\{\varepsilon,w_{e_i}\}\bigr\}
    =
    \chi_{e_i}.
\end{equation}
Finally, $-\Delta^-_{e_k}\le 0=\chi_{e_k}$. Substituting these bounds into the regrouped objective yields
\[
    f^\top\bigl(w+\Delta^+-\Delta^-\bigr)
    \le
    \sum_{i=1}^k w_{e_i} + c_s + \sum_{i=1}^k \chi_{e_i}
    =
    (w^{\mathrm{wc}})^\top f + c_s.
\]
Hence
\begin{equation}
\label{eq:pathwise-upper-infty-short}
    \max_{(\Delta^+,\Delta^-)\in \mathcal D^{\mathrm{S},\infty}(\varepsilon)}
    f^\top\bigl(w+\Delta^+-\Delta^-\bigr)
    \le
    (w^{\mathrm{wc}})^\top f + c_s.
\end{equation}

\smallskip
\noindent\emph{Lower bound.}
To prove the matching lower bound, define $\rho_i\coloneqq \chi_{e_i}$ for $i=1,\dots,k-1$. By definition of $\chi_{e_i}$,
\[
    0\le \rho_i\le \sum_{e\in E_{\mathrm{in}}(v_i)\setminus\{e_i\}} \min\{\varepsilon,w_e\},
\]
and similarly
\[
    0\le c_s\le \sum_{e\in E_{\mathrm{in}}(s)} \min\{\varepsilon,w_e\}.
\]
Hence, for each internal node $v_i$, one can choose nonnegative values on the edges in $ E_{\mathrm{in}}(v_i)\setminus\{e_i\}$, each bounded by $\min\{\varepsilon,w_e\}$, whose sum is exactly $\rho_i$; likewise, one can choose nonnegative values on the edges in $ E_{\mathrm{in}}(s)$, each bounded by $\min\{\varepsilon,w_e\}$, whose sum is exactly $c_s$. Using these choices, define $\Delta^{-,\star}$ by placing exactly those amounts on the corresponding incoming off-path edges at each internal node and on the incoming edges of $s$, and set all remaining components of $\Delta^{-,\star}$ to zero. Define $\Delta^{+,\star}$ by setting
\[
    \Delta^{+,\star}_{e_1}=c_s,
    \qquad
    \Delta^{+,\star}_{e_{i+1}}=\rho_i=\chi_{e_i}\quad (i=1,\dots,k-1),
\]
and all remaining components to zero. By construction, all components are nonnegative, $\Delta^{-,\star}\le w$, and $\|(\Delta^{+,\star},\Delta^{-,\star})\|_\infty\le \varepsilon$. Conservation holds at $s$ and at each internal node because the total assigned inflow equals the total assigned outflow there, and it holds trivially elsewhere. Thus
\[
    (\Delta^{+,\star},\Delta^{-,\star})\in \mathcal D^{\mathrm{S},\infty}(\varepsilon).
\]
Since $\Delta^{-,\star}$ vanishes on all path edges and the only nonzero inflows on path edges are $c_s$ on $e_1$ and $\chi_{e_i}$ on $e_{i+1}$ for $i=1,\dots,k-1$, we obtain
\[
\begin{aligned}
    f^\top\bigl(w+\Delta^{+,\star}-\Delta^{-,\star}\bigr)
    &=
    \sum_{i=1}^k w_{e_i} + c_s + \sum_{i=1}^{k-1}\chi_{e_i} \\
    &=
    \sum_{i=1}^k \bigl(w_{e_i}+\chi_{e_i}\bigr)+c_s \\
    &=
    (w^{\mathrm{wc}})^\top f + c_s,
\end{aligned}
\]
since $\chi_{e_k}=0$. Therefore,
\[
    \max_{(\Delta^+,\Delta^-)\in \mathcal D^{\mathrm{S},\infty}(\varepsilon)}
    f^\top\bigl(w+\Delta^+-\Delta^-\bigr)
    \ge
    (w^{\mathrm{wc}})^\top f + c_s.
\]
Together with \eqref{eq:pathwise-upper-infty-short}, this proves \eqref{eq:pathwise-identity-infty}.

Equation~\eqref{eq:pathwise-identity-infty} is the key pathwise equivalence: for every fixed $s$--$t$ path, its robust value under short-term local diffusion is exactly its value in the precomputed worst-case graph, up to the path-independent source correction term $c_s$.

\medskip
\noindent\emph{Step 2: Optimality of the algorithm.}
Since $c_s$ is independent of the path,
\[
    \min_{f\in\mathcal P_{s,t}}
    \max_{(\Delta^+,\Delta^-)\in \mathcal D^{\mathrm{S},\infty}(\varepsilon)}
    f^\top\bigl(w+\Delta^+-\Delta^-\bigr)
    =
    \min_{f\in\mathcal P_{s,t}} \bigl((w^{\mathrm{wc}})^\top f + c_s\bigr)
    =
    c_s + \min_{f\in\mathcal P_{s,t}} (w^{\mathrm{wc}})^\top f.
\]
Thus, minimizing the robust objective is exactly the same as computing a shortest $s$--$t$ path in the graph with edge weights $w^{\mathrm{wc}}$. Algorithm~\ref{alg:diff-rsp-short-term} does exactly this. Hence it returns an optimal path-incidence vector $f^\star$ together with the optimal robust value.

\medskip
\noindent\emph{Step 3: Running time.}
The quantities $T_u$ are computed by one pass over the incoming adjacency lists, which takes $O(| E|+| V|)$ time. Once these are available, each surcharge $\chi_e$ and each modified weight $w_e^{\mathrm{wc}}$ is computed in constant time per edge, so the total preprocessing cost is $O(| E|+| V|)$. Since $w^{\mathrm{wc}}\ge 0$, the final shortest-path computation can be performed using Dijkstra's algorithm in $O(| E|+| V|\log| V|)$ time. Thus the overall running time is
\[
    O(| E|+| V|\log| V|).
\]
This proves the claim under $\mathcal D^{\mathrm{S},\infty}(\varepsilon)$.

\bigskip
\noindent\textbf{(ii) Case $\mathcal D^{\mathrm{S},1}(\varepsilon)$.} As above, we first establish correctness, and then derive the running-time bound.

\medskip
\noindent\emph{Step 1: Pathwise equivalence.}
Fix any $f\in\mathcal P_{s,t}$, and let the corresponding directed $s$--$t$ path be
\[
    s=v_0 \xrightarrow{e_1} v_1 \xrightarrow{e_2} \cdots \xrightarrow{e_k} v_k=t.
\]
We claim that
\begin{equation}
\label{eq:pathwise-identity-l1}
    \max_{(\Delta^+,\Delta^-)\in \mathcal D^{\mathrm{S},1}(\varepsilon)}
    f^\top\bigl(w+\Delta^+-\Delta^-\bigr)
    =
    \min\Bigl\{w^\top f+\frac{\varepsilon}{2},\; (w^{\mathrm{wc}})^\top f + T_s\Bigr\}.
\end{equation}
We again prove \eqref{eq:pathwise-identity-l1} by matching upper and lower bounds.

\smallskip
\noindent\emph{Upper bound.}
Fix any feasible $(\Delta^+,\Delta^-)\in\mathcal D^{\mathrm{S},1}(\varepsilon)$. Since $\|(\Delta^+,\Delta^-)\|_1\le \varepsilon$ and all components are nonnegative, conservation implies
\[
    \sum_{e\in E}\Delta^+_e
    =
    \sum_{e\in E}\Delta^-_e
    \le
    \frac{\varepsilon}{2}.
\]
Therefore,
\[
    f^\top\bigl(w+\Delta^+-\Delta^-\bigr)
    =
    w^\top f + \sum_{e\in E} f_e\Delta^+_e - \sum_{e\in E} f_e\Delta^-_e
    \le
    w^\top f + \sum_{e\in E} \Delta^+_e
    \le
    w^\top f + \frac{\varepsilon}{2}.
\]
This gives the first upper bound. For the second upper bound, the regrouping identity \eqref{eq:rsp-path-regroup} from the $\mathcal D^{\mathrm{S},\infty}(\varepsilon)$ case holds verbatim:
\[
    f^\top\bigl(w+\Delta^+-\Delta^-\bigr)
    =
    \sum_{i=1}^k w_{e_i}
    + \Delta^+_{e_1}
    + \sum_{i=1}^{k-1}\bigl(\Delta^+_{e_{i+1}}-\Delta^-_{e_i}\bigr)
    - \Delta^-_{e_k}.
\]
At the source, exactly as before,
\[
    \Delta^+_{e_1}
    \le
    \sum_{e\in E_{\mathrm{in}}(s)} \Delta^-_e
    \le
    \sum_{e\in E_{\mathrm{in}}(s)} \min\{\varepsilon,w_e\}
    =
    T_s.
\]
Moreover, for each internal node $v_i$, the same argument as in \eqref{eq:rsp-node-bound} gives $\Delta^+_{e_{i+1}}-\Delta^-_{e_i}\le \chi_{e_i}$, for $i=1,\dots,k-1$. Finally, $-\Delta^-_{e_k}\le 0=\chi_{e_k}$. Substituting these bounds yields
\[
    f^\top\bigl(w+\Delta^+-\Delta^-\bigr)
    \le
    \sum_{i=1}^k w_{e_i} + T_s + \sum_{i=1}^k \chi_{e_i}
    =
    (w^{\mathrm{wc}})^\top f + T_s.
\]
Combining the two upper bounds, we obtain
\begin{equation}
\label{eq:pathwise-upper-l1}
    \max_{(\Delta^+,\Delta^-)\in \mathcal D^{\mathrm{S},1}(\varepsilon)}
    f^\top\bigl(w+\Delta^+-\Delta^-\bigr)
    \le
    \min\Bigl\{w^\top f+\frac{\varepsilon}{2},\; (w^{\mathrm{wc}})^\top f + T_s\Bigr\}.
\end{equation}

\smallskip
\noindent\emph{Matching feasible construction.}
Set
\[
    \Gamma
    \;\coloneqq\;
    \min\Bigl\{\frac{\varepsilon}{2},\; T_s+\sum_{i=1}^k \chi_{e_i}\Bigr\}.
\]
We construct a feasible diffusion attaining value $w^\top f+\Gamma$, which will match the right-hand side of \eqref{eq:pathwise-upper-l1}. Choose nonnegative numbers $\gamma_0,\gamma_1,\dots,\gamma_{k-1}$ such that
\[
    0\le \gamma_0\le T_s,\qquad
    0\le \gamma_i\le \chi_{e_i}\ \ (i=1,\dots,k-1),
    \qquad\text{and}\qquad
    \gamma_0+\sum_{i=1}^{k-1}\gamma_i=\Gamma.
\]
Such a choice exists because
\[
    \Gamma\le T_s+\sum_{i=1}^k \chi_{e_i}
    =
    T_s+\sum_{i=1}^{k-1}\chi_{e_i},
\]
using $\chi_{e_k}=0$. As in the $\mathcal D^{\mathrm{S},\infty}(\varepsilon)$ case, for the source $s$ choose nonnegative values on edges in $ E_{\mathrm{in}}(s)$, each bounded by $\min\{\varepsilon,w_e\}$, whose sum is $\gamma_0$. Likewise, for each internal node $v_i$, choose nonnegative values on edges in $ E_{\mathrm{in}}(v_i)\setminus\{e_i\}$, each bounded by $\min\{\varepsilon,w_e\}$, whose sum is $\gamma_i$. Such choices are possible because $\gamma_0\le T_s$ and $\gamma_i\le \chi_{e_i}\le \sum_{e\in E_{\mathrm{in}}(v_i)\setminus\{e_i\}}\min\{\varepsilon,w_e\}$.

Using these choices, define $\Delta^{-,\star}$ by placing exactly those amounts on the corresponding incoming off-path edges at the source and at each internal node, and set all remaining components of $\Delta^{-,\star}$ to zero. Define $\Delta^{+,\star}$ by setting
\[
    \Delta^{+,\star}_{e_1}=\gamma_0,
    \qquad
    \Delta^{+,\star}_{e_{i+1}}=\gamma_i\quad (i=1,\dots,k-1),
\]
and all remaining components to zero. By construction, all components are nonnegative. Every nonzero component of $\Delta^{-,\star}$ is at most $\min\{\varepsilon,w_e\}$, so $\Delta^{-,\star}\le w$. Conservation holds at the source and at each internal node because, by construction, the total assigned inflow equals the total assigned outflow there, and it holds trivially elsewhere. Finally,
\[
    \sum_{e\in E}\Delta^{+,\star}_e
    =
    \gamma_0+\sum_{i=1}^{k-1}\gamma_i
    =
    \Gamma,
\]
and likewise
\[
    \sum_{e\in E}\Delta^{-,\star}_e=\Gamma.
\]
Hence
\[
    \|(\Delta^{+,\star},\Delta^{-,\star})\|_1
    =
    \sum_{e\in E}\Delta^{+,\star}_e+\sum_{e\in E}\Delta^{-,\star}_e
    =
    2\Gamma
    \le
    \varepsilon,
\]
so $(\Delta^{+,\star},\Delta^{-,\star})\in\mathcal D^{\mathrm{S},1}(\varepsilon)$. Moreover, $\Delta^{-,\star}$ vanishes on all path edges, while the total inflow placed on path edges is exactly
\[
    \gamma_0+\sum_{i=1}^{k-1}\gamma_i=\Gamma.
\]
Therefore,
\[
    f^\top\bigl(w+\Delta^{+,\star}-\Delta^{-,\star}\bigr)
    =
    w^\top f + \Gamma
    =
    w^\top f + \min\Bigl\{\frac{\varepsilon}{2},\; T_s+\sum_{i=1}^k \chi_{e_i}\Bigr\} 
    =
    \min\Bigl\{w^\top f+\frac{\varepsilon}{2},\; (w^{\mathrm{wc}})^\top f + T_s\Bigr\}.
\]
This matches the upper bound in \eqref{eq:pathwise-upper-l1}, and proves \eqref{eq:pathwise-identity-l1}.

Equation~\eqref{eq:pathwise-identity-l1} is the corresponding pathwise equivalence for the global-budget case: for every fixed $s$--$t$ path, its robust value under short-term $\ell_1$ diffusion is the smaller of the budget-limited value and the value induced by the precomputed worst-case graph.

\medskip
\noindent\emph{Step 2: Optimality of the algorithm.}
By \eqref{eq:pathwise-identity-l1},
\[
\begin{aligned}
    \min_{f\in\mathcal P_{s,t}}
    \max_{(\Delta^+,\Delta^-)\in \mathcal D^{\mathrm{S},1}(\varepsilon)}
    f^\top\bigl(w+\Delta^+-\Delta^-\bigr)
    &=
    \min_{f\in\mathcal P_{s,t}}
    \min\Bigl\{w^\top f+\frac{\varepsilon}{2},\; (w^{\mathrm{wc}})^\top f + T_s\Bigr\}.
\end{aligned}
\]
Since, for any two real-valued functions $A$ and $B$, one has
\[
    \min_f \min\{A(f),B(f)\}
    =
    \min\Bigl\{\min_f A(f),\min_f B(f)\Bigr\},
\]
the robust optimum is obtained by comparing
\[
    \min_{f\in\mathcal P_{s,t}} w^\top f + \frac{\varepsilon}{2}
    \qquad\text{and}\qquad
    \min_{f\in\mathcal P_{s,t}} (w^{\mathrm{wc}})^\top f + T_s.
\]
Algorithm~\ref{alg:diff-rsp-short-term} computes exactly these two candidate values and returns the smaller one. Hence it returns an optimal path-incidence vector $f^\star$ together with the optimal robust value under $\mathcal D^{\mathrm{S},1}(\varepsilon)$.

\medskip
\noindent\emph{Step 3: Running time.}
The preprocessing cost is $O(| E|+| V|)$, exactly as in the $\mathcal D^{\mathrm{S},\infty}(\varepsilon)$ case. The algorithm then solves two shortest-path instances, one with edge weights $w^{\mathrm{wc}}$ and one with edge weights $w$. Thus the overall running time remains
\[
    O(| E|+| V|\log| V|).
\]
This proves the claim under $\mathcal D^{\mathrm{S},1}(\varepsilon)$, and completes the proof of the proposition.

\subsection{Proof of Proposition~\ref{prop:rsp-linf-hardness}}

Fix an instance of MinSAT with formula $\Phi=C_1\wedge\cdots\wedge C_m$ over variables $x_1,\ldots,x_n$, and let $G'=(V',E')$ be the graph constructed using the rules described in Section~\ref{subsec:rsp-np-local} and Figure~\ref{fig:rsp-linf-construction}. We may assume without loss of generality that every clause is nonempty, since an empty clause is never satisfied and can be removed without changing $\mathrm{sat}(A)$ for any assignment $A$. Recall that the budget is $\varepsilon=1$, that the only edges of positive nominal weight are the clause source edges $\alpha_j=(d_j,c_j)$ and the blocker source edges $\eta_i=(\bar d_i,\bar c_i)$, each with nominal weight $1$, and that every other edge in $G'$ has nominal weight $0$. We prove the claim in four steps.

\medskip
\noindent\emph{Step 1: feasible $s$--$t$ paths in $G'$ encode truth assignments.}
We first show that every directed $s$--$t$ path in $G'$ remains on the variable backbone and chooses exactly one branch from each pair $(T_i,F_i)$. Indeed, the dummy vertices $d_j$ and $\bar d_i$ have no incoming edges. Moreover, the only incoming edge of $c_j$ is the clause source edge $\alpha_j=(d_j,c_j)$, and the only incoming edge of $\bar c_i$ is the blocker source edge $\eta_i=(\bar d_i,\bar c_i)$. Since no directed path starting at $s$ can reach any dummy vertex, it cannot reach any vertex $c_j$ or $\bar c_i$. Consequently, no directed $s$--$t$ path can use a source edge, nor any connector edge leaving a vertex $c_j$ or $\bar c_i$.

It follows that any directed $s$--$t$ path must pass through the variable gadgets in their serial order. Within gadget $i$, the only directed routes from $a_i$ to $b_i$ are the two branches $T_i$ and $F_i$. After reaching $b_i$, the path must use the blocking edge $h_i$, which leads to $a_{i+1}$ for $i<n$ and to $t$ for $i=n$. Hence every directed $s$--$t$ path $P$ chooses exactly one branch from each pair $(T_i,F_i)$.

We associate with $P$ the truth assignment $A(P)$ defined by
\[
    A(P)_i =
    \begin{cases}
        \mathrm{true}, & \text{if } P \text{ uses } T_i,\\
        \mathrm{false}, & \text{if } P \text{ uses } F_i.
    \end{cases}
\]
Conversely, every truth assignment $A$ defines a canonical path $P(A)$ by choosing $T_i$ when $A_i=\mathrm{true}$ and $F_i$ when $A_i=\mathrm{false}$, for each $i=1,\ldots,n$. Finally, no directed $s$--$t$ path uses a source edge, and source edges are the only edges of positive nominal weight. Thus, for every directed $s$--$t$ path $P$ with incidence vector $f$, we have
\begin{equation}
\label{eq:rsp-linf-nominal-zero-final}
    f^\top w=0.
\end{equation}

\medskip
\noindent\emph{Step 2: lower bound on the robust value.}
Fix a directed $s$--$t$ path $P$ in $G'$ and let $A=A(P)$ be the assignment encoded by $P$. We construct a feasible diffusion $(\Delta^{+,\star},\Delta^{-,\star})\in\mathcal D^{\mathrm{L},\infty}(1)$ such that
\begin{equation}
\label{eq:rsp-linf-lower-final}
    f^\top\bigl(w+\Delta^{+,\star}-\Delta^{-,\star}\bigr)
    \ge n+\mathrm{sat}(A).
\end{equation}

\smallskip
\noindent\emph{(a) Baseline contributions from blocker source edges.}
For each variable gadget $i$, let $\gamma_i$ denote the connector edge from $\bar c_i$ to the tail of the blocking edge $h_i$. Set
\[
    \Delta^{-,\star}_{\eta_i}=1,
    \qquad
    \Delta^{+,\star}_{\gamma_i}=1,
    \qquad
    \Delta^{-,\star}_{\gamma_i}=1,
    \qquad
    \Delta^{+,\star}_{h_i}=1.
\]
This routes one unit from the blocker source edge $\eta_i$ to the blocking edge $h_i$.

\smallskip
\noindent\emph{(b) Clause contributions from satisfied clauses.}
For each clause $C_j$ satisfied by $A$, choose one literal occurrence in $C_j$ that is true under $A$. The corresponding literal-occurrence edge lies on the path $P$; denote this edge by $e_j$. Let $\kappa_j$ denote the connector edge from $c_j$ to the tail of $e_j$. Set
\[
    \Delta^{-,\star}_{\alpha_j}=1,
    \qquad
    \Delta^{+,\star}_{\kappa_j}=1,
    \qquad
    \Delta^{-,\star}_{\kappa_j}=1,
    \qquad
    \Delta^{+,\star}_{e_j}=1.
\]
For clauses not satisfied by $A$, all components associated with their clause source edges are set to zero. All remaining components of $\Delta^{+,\star}$ and $\Delta^{-,\star}$ are also set to zero.

We verify feasibility of this construction:

\smallskip
\noindent\emph{Budget constraint.}
The local $\ell_\infty$ budget holds because every component of $(\Delta^{+,\star},\Delta^{-,\star})$ is either $0$ or $1$. 

\smallskip
\noindent\emph{Propagation constraint.}
For each used source edge, namely $e=\eta_i$ for some $i\in\{1,\ldots,n\}$ or $e=\alpha_j$ for some satisfied clause $C_j$, we have $\Delta^{-,\star}_e=1$, $w_e=1$, and $\Delta^{+,\star}_e=0$, and hence
\[
    \Delta^{-,\star}_e
    =
    1
    =
    w_e
    \le
    w_e+\Delta^{+,\star}_e.
\]
For each used connector edge, namely $e=\gamma_i$ for some $i\in\{1,\ldots,n\}$ or $e=\kappa_j$ for some satisfied clause $C_j$, we have $\Delta^{-,\star}_e=\Delta^{+,\star}_e=1$ and $w_e=0$, and hence
\[
    \Delta^{-,\star}_e
    =
    1
    =
    \Delta^{+,\star}_e
    =
    w_e+\Delta^{+,\star}_e.
\]
Finally, for each receiving path edge, namely $e=h_i$ for some $i\in\{1,\ldots,n\}$ or $e=e_j$ for some satisfied clause $C_j$, we have $\Delta^{-,\star}_e=0$, and therefore
\[
    \Delta^{-,\star}_e
    =
    0
    \le
    w_e+\Delta^{+,\star}_e.
\]
All unused edges have $\Delta^{+,\star}_e=\Delta^{-,\star}_e=0$, so the propagation constraint is satisfied on every edge.

\smallskip
\noindent\emph{Conservation.} 
At each used clause-source vertex $c_j$, one unit enters through $\alpha_j$ and one unit leaves through the selected connector $\kappa_j$. At each blocker-source vertex $\bar c_i$, one unit enters through $\eta_i$ and one unit leaves through $\gamma_i$. At the tail of each selected literal-occurrence edge $e_j$, one unit enters through $\kappa_j$ and one unit leaves through $e_j$. Similarly, at the tail of each blocking edge $h_i$, one unit enters through $\gamma_i$ and one unit leaves through $h_i$. All other nodes carry zero perturbation flow. Each constructed unit follows a source edge, then a connector edge, and finally enters one receiving path edge; the conservation equalities above hold at the two intermediate vertices of each such route. Therefore
\[
    (\Delta^{+,\star},\Delta^{-,\star})\in\mathcal D^{\mathrm{L},\infty}(1).
\]

Finally, the only positive contributions on the path $P$ are the $n$ contributions on the blocking edges $h_1,\ldots,h_n$ and the contributions on the selected literal-occurrence edges $e_j$ for satisfied clauses. Since no path edge has positive $\Delta^{-,\star}$, \eqref{eq:rsp-linf-nominal-zero-final} gives
\[
\begin{aligned}
    f^\top\bigl(w+\Delta^{+,\star}-\Delta^{-,\star}\bigr)
    &=
    f^\top(\Delta^{+,\star}-\Delta^{-,\star}) \\
    &=
    \sum_{i=1}^n \Delta^{+,\star}_{h_i}
    +
    \sum_{j:\, C_j(A)=\mathrm{true}} \Delta^{+,\star}_{e_j} \\
    &=
    n+\mathrm{sat}(A).
\end{aligned}
\]
This proves \eqref{eq:rsp-linf-lower-final}.

\medskip
\noindent\emph{Step 3: upper bound on the robust value.}
We now show that no feasible diffusion can produce value larger than $n+\mathrm{sat}(A)$. Fix any feasible diffusion $(\Delta^+,\Delta^-)\in\mathcal D^{\mathrm{L},\infty}(1)$. By \eqref{eq:rsp-linf-nominal-zero-final},
\[
    f^\top\bigl(w+\Delta^+-\Delta^-\bigr)
    =
    f^\top(\Delta^+-\Delta^-).
\]
Only first entries from outside the path can increase this total value. Indeed, if an amount $\lambda$ is forwarded from one path edge to the next, then it contributes $-\lambda$ on the edge from which it is removed and $+\lambda$ on the next path edge, for zero net change in $f^\top(\Delta^+-\Delta^-)$. Hence $f^\top(\Delta^+-\Delta^-)$ is at most the total perturbation mass that first enters the path $P$ from outside the path.

We bound this first-entry mass by separating the path edges into three classes: blocking edges, literal-occurrence edges of the chosen branches, and branch edges of the chosen branches that are not literal-occurrence edges. By construction, no connector enters a non-occurrence branch edge. Hence such edges cannot be first-entry edges from outside the path. Thus only the first two classes can contribute first-entry mass.

\smallskip
\noindent\emph{Blocking-edge entries.}
The path $P$ contains exactly the blocking edges $h_1,\ldots,h_n$. Since the local budget is $\varepsilon=1$, each satisfies $\Delta^+_{h_i}\le 1$. Therefore the total first-entry mass into blocking edges is at most
\[
    \sum_{i=1}^n \Delta^+_{h_i}\le n.
\]

\smallskip
\noindent\emph{Literal-occurrence entries.}
Consider first-entry mass into literal-occurrence edges of the chosen branches. Such mass can only originate from clause source edges. Indeed, blocker source edges have connectors only to the blocking edges $h_i$, whereas the connectors leaving $c_j$ enter exactly the literal-occurrence edges corresponding to literals contained in $C_j$.

Each clause source edge can emit at most one unit in total. To see this, fix $j$. Since $d_j$ has no incoming edges, conservation at $d_j$ gives $\Delta^+_{\alpha_j}=0$. The propagation constraint on $\alpha_j$ then gives
\[
    \Delta^-_{\alpha_j}\le w_{\alpha_j}+\Delta^+_{\alpha_j}=1.
\]
Since the only incoming edge of $c_j$ is $\alpha_j$, conservation at $c_j$ implies that the total amount sent from $c_j$ through all clause connectors is at most $\Delta^-_{\alpha_j}\le 1$.

Now suppose that $C_j$ is not satisfied by the assignment $A=A(P)$. Then every literal in $C_j$ is false under $A$. Hence every connector leaving $c_j$ enters a literal-occurrence edge on an unchosen branch. Since the tails of literal-occurrence edges are private internal branch vertices, source $\alpha_j$ cannot first enter $P$ through a literal-occurrence edge of a chosen branch. If any mass from this source reaches $P$ later, it must do so after traveling along an unchosen branch and first entering the path through a blocking edge $h_i$ for some $i$. But first entries through blocking edges have already been bounded by the baseline contribution $\sum_{i=1}^n \Delta^+_{h_i}\le n$. Thus an unsatisfied clause source can be absorbed by the blocker baseline, but it cannot create an additional literal-occurrence contribution. By contrast, if $C_j$ is satisfied by $A$, then source $\alpha_j$ can contribute at most one unit of first-entry mass into the literal-occurrence edges of $P$. Therefore the total first-entry mass into literal-occurrence edges of the chosen branches is at most $\mathrm{sat}(A)$.

Combining the two bounds, we obtain
\[
    f^\top(\Delta^+-\Delta^-)
    \le
    n+\mathrm{sat}(A).
\]
Since $(\Delta^+,\Delta^-)$ was arbitrary, this proves
\begin{equation}
\label{eq:rsp-linf-upper-final}
    \max_{(\Delta^+,\Delta^-)\in\mathcal D^{\mathrm{L},\infty}(1)}
    f^\top\bigl(w+\Delta^+-\Delta^-\bigr)
    \le
    n+\mathrm{sat}(A).
\end{equation}

Combining \eqref{eq:rsp-linf-lower-final} and \eqref{eq:rsp-linf-upper-final}, we conclude that, for every directed $s$--$t$ path $P$ in $G'$,
\[
    \max_{(\Delta^+,\Delta^-)\in\mathcal D^{\mathrm{L},\infty}(1)}
    f^\top\bigl(w+\Delta^+-\Delta^-\bigr)
    =
    n+\mathrm{sat}(A(P)).
\]

\medskip
\noindent\emph{Step 4: conclusion of the reduction.}
The preceding equality implies that, for every integer $K\ge 0$,
\[
    \exists \text{ a truth assignment } A \text{ with } \mathrm{sat}(A)\le K
\]
if and only if
\[
    \exists \text{ a directed $s$--$t$ path } P \text{ in } G'
    \text{ with robust value } \le n+K.
\]
Indeed, the forward implication uses the canonical path $P(A)$, while the reverse implication uses the assignment $A(P)$ induced by any feasible path $P$. Thus the constructed Diff-RSP instance recovers the threshold structure of MinSAT.

It remains only to verify the running time of the construction. The variable backbone contributes $O(n)$ vertices and edges. The literal branches contribute $O(1)$ vertices and edges per literal occurrence, together with $O(1)$ auxiliary vertices and edges per variable branch, and therefore have total size $O(n+\ell_\Phi)$. The clause source edges contribute $m$ additional source edges and one connector for each literal occurrence, for total size $O(m+\ell_\Phi)$. Since every clause is nonempty, $m\le \ell_\Phi$. The blocker source edges contribute $O(n)$ additional vertices and edges. Hence the graph $G'$ has size $O(n+\ell_\Phi)$ and can be constructed from $\Phi$ in $O(n+\ell_\Phi)$ time. This completes the proof.

\subsection{Proof of Proposition~\ref{prop:rsp-l1-hardness}}

Fix an instance of Most Secluded Path on a directed graph $G=(V,E)$ with terminals $s,t\in V$, and let $G''=(V'',E'')$ be the graph constructed using the rules described in Section~\ref{subsec:rsp-np-global} and Figure~\ref{fig:rsp-construction-Gdoubleprime}. Recall that the budget is $\varepsilon=4|V|$, that the only edges of positive nominal weight in $G''$ are the dummy edges $\alpha_v=(d_v,c_v)$, each with nominal weight $1$, and that every other edge in $G''$ has nominal weight $0$. We prove the claim in four steps.

\medskip
\noindent\emph{Step 1: feasible $s'$--$t'$ paths in $G''$ are canonical.}
Set $s' \coloneqq s_1$ and $t' \coloneqq t_{4|V|+1}$. We first show that every directed $s'$--$t'$ path in $G''$ is of the form
\[
    E_{q_1}\to a_{q_1,q_2}\to E_{q_2}\to \cdots \to a_{q_{k-1},q_k}\to E_{q_k},
\]
for some directed $s$--$t$ path $Q=(q_1=s,q_2,\ldots,q_k=t)$ in $G$. We call such a path \emph{canonical}. Indeed, the only incoming edge of the entry node $c_u$ is the dummy edge $\alpha_u=(d_u,c_u)$, and the dummy vertex $d_u$ has no incoming edges. Since $s'=s_1$ lies in the chain gadget $E_s$, no directed path starting at $s'$ can ever reach a node $c_u$. Consequently, no directed $s'$--$t'$ path can use any connector edge, because every connector edge leaves from some $c_v$. It follows that the only edges that can move a directed path from one gadget to another are the anchor edges.

By construction, an anchor edge $a_{u,v}$ exists if and only if $(u,v)\in E$, and it leaves the terminal vertex $u_{4|V|+1}$ of the chain gadget $E_u$ and enters the initial vertex $v_1$ of the chain gadget $E_v$. Moreover, once a directed path enters a gadget at $v_1$, the only outgoing edge from $v_i$ is the next chain edge
\[
    e_v^i=(v_i,v_{i+1}),
    \qquad i=1,\ldots,4|V|.
\]
Hence the path must traverse the entire gadget $E_v$ before it can leave through an anchor edge. Therefore, every directed $s'$--$t'$ path in $G''$ traverses a sequence of full chain gadgets linked by anchor edges, and this sequence is exactly a directed $s$--$t$ path in $G$.

Now fix an arbitrary directed $s$--$t$ path
\[
    Q=(q_1=s,q_2,\ldots,q_k=t)
\]
in $G$, and let $P(Q)$ denote the corresponding canonical directed $s'$--$t'$ path in $G''$. Let $f$ be its incidence vector. Since no canonical path uses any dummy edge, and dummy edges are the only edges of positive nominal weight, we immediately have
\begin{equation}
\label{eq:rsp-l1-nominal-zero-final}
    f^\top w = 0.
\end{equation}

\medskip
\noindent\emph{Step 2: lower bound on the robust value.}
We now construct a feasible diffusion $(\Delta^{+,\star},\Delta^{-,\star})\in \mathcal D^{\mathrm{L},1}(4|V|)$ such that
\begin{equation}
\label{eq:rsp-l1-lower-final}
    f^\top\bigl(w+\Delta^{+,\star}-\Delta^{-,\star}\bigr)
    \ge |N[V(Q)]|.
\end{equation}
For each node $u\in N[V(Q)]\setminus V(Q)$, choose one assigned receiver
\[
    r(u)\in V(Q)
    \qquad\text{such that}\qquad
    u\in N[r(u)]\setminus\{r(u)\}.
\]
Such a choice is possible by the definition of the closed out-neighborhood. We fix one such choice for every off-path exposed node. Thus each node in $N[V(Q)]\setminus V(Q)$ is assigned to exactly one on-path receiver, even if it belongs to $N[v]\setminus\{v\}$ for more than one node $v\in V(Q)$. Consequently, the sets
\[
    \{u\in N[V(Q)]\setminus V(Q): r(u)=v\},
    \qquad v\in V(Q),
\]
are pairwise disjoint and their union is $N[V(Q)]\setminus V(Q)$. For each on-path node $v\in V(Q)$, let
\[
    m_v \coloneqq
    \bigl|\{u\in N[V(Q)]\setminus V(Q): r(u)=v\}\bigr|
\]
denote the number of exposed off-path nodes assigned to $v$.

We define $(\Delta^{+,\star},\Delta^{-,\star})$ by specifying its nonzero components:

\smallskip
\noindent\emph{(a) Self-contributions of on-path nodes.}
For each $v\in V(Q)$, set
\[
    \Delta^{-,\star}_{\alpha_v}=1,
    \qquad
    \Delta^{+,\star}_{\beta_v}=1,
    \qquad
    \Delta^{-,\star}_{\beta_v}=1,
    \qquad
    \Delta^{+,\star}_{e_v^1}=1.
\]

\smallskip
\noindent\emph{(b) Contributions of off-path exposed nodes.}
Fix $u\in N[V(Q)]\setminus V(Q)$, and let $v=r(u)$. Since $u\in N[v]\setminus\{v\}$, the corresponding connector edge in the construction is $c_{u,v}\coloneqq (c_u,v_2)$. Set
\[
    \Delta^{-,\star}_{\alpha_u}=1,
    \qquad
    \Delta^{+,\star}_{c_{u,v}}=1,
    \qquad
    \Delta^{-,\star}_{c_{u,v}}=1.
\]
In addition, for each $v\in V(Q)$, set
\[
    \Delta^{+,\star}_{e_v^2}=m_v.
\]
All remaining components of $\Delta^{+,\star}$ and $\Delta^{-,\star}$ are set to zero.

We verify the feasibility of this construction:

\smallskip
\noindent\emph{Budget constraint.}
Each on-path node $v\in V(Q)$ contributes exactly
\[
    \Delta^{-,\star}_{\alpha_v}
    +
    \Delta^{+,\star}_{\beta_v}
    +
    \Delta^{-,\star}_{\beta_v}
    +
    \Delta^{+,\star}_{e_v^1}
    =4
\]
to the $\ell_1$ norm. The off-path exposed contributions have total budget
\[
\begin{aligned}
    \sum_{u\in N[V(Q)]\setminus V(Q)}
    \left(
        \Delta^{-,\star}_{\alpha_u}
        +
        \Delta^{+,\star}_{c_{u,r(u)}}
        +
        \Delta^{-,\star}_{c_{u,r(u)}}
    \right)
    +
    \sum_{v\in V(Q)}\Delta^{+,\star}_{e_v^2}  
    &=
    3\bigl(|N[V(Q)]|-|V(Q)|\bigr)
    +
    \sum_{v\in V(Q)}m_v \\
    &=
    4\bigl(|N[V(Q)]|-|V(Q)|\bigr).
\end{aligned}
\]
Together with the $4|V(Q)|$ units used by the self-contributions of the on-path nodes, this gives
\[
    \|(\Delta^{+,\star},\Delta^{-,\star})\|_1
    =
    4|N[V(Q)]|
    \le 4|V|
    =
    \varepsilon.
\]

\smallskip
\noindent\emph{Propagation constraint.}
On each dummy edge $\alpha_u$, we have
\[
    \Delta^{-,\star}_{\alpha_u}=1=w_{\alpha_u}\le w_{\alpha_u}+\Delta^{+,\star}_{\alpha_u}.
\]
On each used entry edge $\beta_v$ and each used connector edge $c_{u,v}$, we have $\Delta^{-,\star}_e = 1 = \Delta^{+,\star}_e$, and these edges have nominal weight $0$, so
\[
    \Delta^{-,\star}_e \le w_e+\Delta^{+,\star}_e.
\]
Finally, on the used path edges $e_v^1$ and $e_v^2$, we have
\[
    \Delta^{-,\star}_{e_v^1}=0\le \Delta^{+,\star}_{e_v^1},
    \qquad
    \Delta^{-,\star}_{e_v^2}=0\le \Delta^{+,\star}_{e_v^2},
\]
and these edges also have nominal weight $0$.

\smallskip
\noindent\emph{Conservation.}
At each entry node $c_u$, the only incoming edge is $\alpha_u$. If $u\in V(Q)$, then the only nonzero outgoing term is $\Delta^{+,\star}_{\beta_u}=1$. If
\[
    u\in N[V(Q)]\setminus V(Q),
\]
then the only nonzero outgoing term is $\Delta^{+,\star}_{c_{u,r(u)}}=1$. Thus conservation holds at every node $c_u$.

At each initial gadget vertex $v_1$ with $v\in V(Q)$, the only nonzero incoming term is $\Delta^{-,\star}_{\beta_v}=1$, and the only nonzero outgoing term is $\Delta^{+,\star}_{e_v^1}=1$. Hence conservation holds at $v_1$. At each second gadget vertex $v_2$ with $v\in V(Q)$, the nonzero incoming terms are exactly
\[
    \Delta^{-,\star}_{c_{u,v}}=1
    \qquad
    \text{for all }u\in N[V(Q)]\setminus V(Q)\text{ with }r(u)=v,
\]
whose total is $m_v$, and the only nonzero outgoing term is
\[
    \Delta^{+,\star}_{e_v^2}=m_v.
\]
Hence conservation also holds at $v_2$. All remaining nodes carry zero flow and satisfy conservation trivially. Therefore,
\[
    (\Delta^{+,\star},\Delta^{-,\star})\in \mathcal D^{\mathrm{L},1}(4|V|).
\]

Finally, the only nonzero terms on path edges are the $|V(Q)|$ self-contributions on the edges $e_v^1$ and the $|N[V(Q)]|-|V(Q)|$ exposed off-path contributions aggregated on the edges $e_v^2$. Therefore
\[
\begin{aligned}
    f^\top\bigl(w+\Delta^{+,\star}-\Delta^{-,\star}\bigr)
    =
    f^\top(\Delta^{+,\star}-\Delta^{-,\star})
    &=
    \sum_{v\in V(Q)} \Delta^{+,\star}_{e_v^1}
    +
    \sum_{v\in V(Q)} \Delta^{+,\star}_{e_v^2} \\
    &=
    |V(Q)|+\sum_{v\in V(Q)} m_v \\
    &=
    |N[V(Q)]|.
\end{aligned}
\]
This proves \eqref{eq:rsp-l1-lower-final}.

\medskip
\noindent\emph{Step 3: upper bound on the robust value.}
We now show that
\begin{equation}
\label{eq:rsp-l1-upper-final}
    \max_{(\Delta^+,\Delta^-)\in \mathcal D^{\mathrm{L},1}(4|V|)}
    f^\top\bigl(w+\Delta^+-\Delta^-\bigr)
    <
    |N[V(Q)]|+\tfrac{1}{2}.
\end{equation}
Fix any feasible diffusion $(\Delta^+,\Delta^-)\in \mathcal D^{\mathrm{L},1}(4|V|)$. Since the canonical path uses no dummy edges, \eqref{eq:rsp-l1-nominal-zero-final} gives
\[
    f^\top\bigl(w+\Delta^+-\Delta^-\bigr)=f^\top(\Delta^+-\Delta^-).
\]
We begin with two observations.

\smallskip
\noindent\emph{Observation 1: each source gadget can emit at most one unit.}
Fix $u\in V$. Since the dummy vertex $d_u$ has no incoming edges, conservation at $d_u$ implies $\Delta^+_{\alpha_u}=0$. Hence the long-term propagation constraint on $\alpha_u$ becomes
\[
    \Delta^-_{\alpha_u}\le w_{\alpha_u}+\Delta^+_{\alpha_u}=1.
\]
Moreover, the only incoming edge of $c_u$ is $\alpha_u$, so conservation at $c_u$ gives
\[
    \Delta^+_{\beta_u}
    +
    \sum_{v:\,u\in N[v]\setminus\{v\}}
    \Delta^+_{c_{u,v}}
    =
    \Delta^-_{\alpha_u}
    \le 1.
\]
Thus the gadget of $u$ can emit at most one unit of perturbation mass in total.

\smallskip
\noindent\emph{Observation 2: after a unit first enters the canonical path, forwarding it further along the path does not increase the net contribution.}
Indeed, if an amount $\lambda$ moves from some path edge $e$ to the next path edge $e'$, then this contributes $-\lambda$ on $e$ and $+\lambda$ on $e'$, for a net change of $0$ in $f^\top(\Delta^+-\Delta^-)$. Therefore, the total value $f^\top(\Delta^+-\Delta^-)$ is at most the total perturbation mass that \emph{first enters} the canonical path $P(Q)$.

We now proceed to prove~\eqref{eq:rsp-l1-upper-final}. For each node $u\in V$, let $\mu_u$ denote the total first-entry mass on the canonical path that originates from the gadget of $u$. By Observation~1,
\[
    0\le \mu_u\le 1,
    \qquad u\in V.
\]
Let
\[
    M_{\mathrm{exp}}
    \coloneqq
    \sum_{u\in N[V(Q)]}\mu_u,
    \qquad
    M_{\mathrm{out}}
    \coloneqq
    \sum_{u\in V\setminus N[V(Q)]}\mu_u.
\]
Then Observation~2 implies
\begin{equation}
\label{eq:rsp-l1-first-entry-final}
    f^\top(\Delta^+-\Delta^-)
    \le
    M_{\mathrm{exp}}+M_{\mathrm{out}}.
\end{equation}
Since there are exactly $|N[V(Q)]|$ exposed nodes and each contributes at most one unit, we also have
\begin{equation}
\label{eq:rsp-l1-exp-upper-final}
    M_{\mathrm{exp}}\le |N[V(Q)]|.
\end{equation}
We next derive a lower bound on the $\ell_1$ budget required to create these first-entry contributions.

\smallskip
\noindent\emph{Exposed sources.}
Fix $u\in N[V(Q)]$, and let $\mu_u=\lambda$. We claim that any perturbation in which $\lambda$ units of mass originating from the gadget of $u$ first enter the canonical path requires at least $4\lambda$ units of $\ell_1$ budget, regardless of the particular routes used. Indeed, the diffusion must use $\lambda$ units of outflow on the dummy edge $\alpha_u$, contributing $\lambda$ to the budget. Before this mass first reaches the canonical path, it must be carried through non-path transfer edges with total flow at least $\lambda$. These transfer edges contribute at least $\lambda$ through $\Delta^+$ and at least $\lambda$ through $\Delta^-$. Finally, the first canonical-path edges reached by this mass receive total inflow $\lambda$, contributing another $\lambda$ to $\Delta^+$.

Hence producing $\lambda$ units of first-entry mass from any exposed source costs at least $\lambda+2\lambda+\lambda=4\lambda$ in $\|(\Delta^+,\Delta^-)\|_1$.

\smallskip
\noindent\emph{Non-exposed sources.}
Fix $u\notin N[V(Q)]$, and let $\mu_u=\lambda$. Since $u$ is not in the exposure set, there is no connector from $c_u$ into any on-path gadget. Indeed, a connector from $c_u$ into the gadget of some $v\in V(Q)$ exists only when $u\in N[v]\setminus\{v\}$, which would imply $u\in N[V(Q)]$. Moreover, once mass leaves $c_u$, it can never reach any node of the form $c_x$, because the only incoming edge of $c_x$ is the dummy edge $\alpha_x$. Therefore, before mass from $u$ can first enter the canonical path, it must first leave $c_u$ through either its own entry edge $\beta_u$ or connector edges into off-path gadgets. 

In either case, each unit of such mass must traverse at least $4|V|-1$ chain edges outside the canonical path before first reaching the canonical path. The bound $4|V|-1$ is the most favorable case for the adversary: it occurs when the mass enters an off-path gadget through a connector, thereby skipping the first chain edge of that gadget. Entering through $\beta_u$ requires traversing the full $4|V|$-edge chain before any anchor can be used. This argument also covers routes that reach the canonical path only after passing through several off-path gadgets; we charge only the first off-path chain segment before first entry into the canonical path.

Thus producing $\lambda$ units of first-entry mass from $u$ requires, at a minimum:
\begin{itemize}
    \item $\lambda$ units of outflow on the dummy edge $\alpha_u$;
    \item non-dummy transfer edges leaving $c_u$ carrying total flow at least $\lambda$, contributing at least $\lambda$ to $\Delta^+$ and at least $\lambda$ to $\Delta^-$;
    \item chain edges outside the canonical path carrying total flow at least $(4|V|-1)\lambda$, contributing at least $(4|V|-1)\lambda$ to $\Delta^+$ and at least $(4|V|-1)\lambda$ to $\Delta^-$;
    \item total inflow $\lambda$ on the first canonical-path edges reached by this mass, contributing $\lambda$ to $\Delta^+$.
\end{itemize}
Hence producing $\lambda$ units of first-entry mass from a non-exposed source costs at least
\[
    \lambda + 2\lambda + 2(4|V|-1)\lambda + \lambda
    =
    (8|V|+2)\lambda
\]
in $\|(\Delta^+,\Delta^-)\|_1$.

Summing over all sources, we obtain the budget lower bound
\begin{equation}
\label{eq:rsp-l1-budget-final}
    4M_{\mathrm{exp}}+(8|V|+2)M_{\mathrm{out}}
    \le
    \|(\Delta^+,\Delta^-)\|_1
    \le
    4|V|.
\end{equation}
Therefore,
\[
    M_{\mathrm{out}}
    \le
    \frac{4|V|-4M_{\mathrm{exp}}}{8|V|+2}.
\]
Combining this with \eqref{eq:rsp-l1-first-entry-final} yields
\[
    f^\top(\Delta^+-\Delta^-)
    \le
    M_{\mathrm{exp}}
    +
    \frac{4|V|-4M_{\mathrm{exp}}}{8|V|+2}.
\]
The right-hand side is increasing in $M_{\mathrm{exp}}$, because
\[
    1-\frac{4}{8|V|+2}>0.
\]
Using \eqref{eq:rsp-l1-exp-upper-final}, we conclude
\[
\begin{aligned}
    f^\top(\Delta^+-\Delta^-)
    &\le
    |N[V(Q)]|
    +
    \frac{4|V|-4|N[V(Q)]|}{8|V|+2}
    <
    |N[V(Q)]|+\frac{1}{2}.
\end{aligned}
\]
This proves \eqref{eq:rsp-l1-upper-final}.

\medskip
\noindent\emph{Step 4: conclusion of the reduction.}
Combining \eqref{eq:rsp-l1-nominal-zero-final}, \eqref{eq:rsp-l1-lower-final}, and \eqref{eq:rsp-l1-upper-final}, we obtain
\[
    |N[V(Q)]|
    \le
    \max_{(\Delta^+,\Delta^-)\in \mathcal D^{\mathrm{L},1}(4|V|)}
    f^\top\bigl(w+\Delta^+-\Delta^-\bigr)
    <
    |N[V(Q)]|+\tfrac{1}{2}
\]
for every directed $s$--$t$ path $Q$ in $G$ and its corresponding canonical path $P(Q)$ in $G''$. Since $|N[V(Q)]|$ is an integer, the preceding bounds imply that, for every integer $K\ge 0$,
\[
    |N[V(Q)]|\le K
    \qquad\Longleftrightarrow\qquad
    \max_{(\Delta^+,\Delta^-)\in \mathcal D^{\mathrm{L},1}(4|V|)}
    f^\top\bigl(w+\Delta^+-\Delta^-\bigr)
    < K+\tfrac{1}{2}
\]
for the canonical path $P(Q)$ corresponding to $Q$. By Step~1, every feasible directed $s'$--$t'$ path in $G''$ is canonical. Therefore, for every integer $K\ge 0$,
\[
    \exists \text{ an $s$--$t$ path } Q \text{ in } G \text{ with } |N[V(Q)]|\le K
\]
if and only if
\[
    \exists \text{ an $s'$--$t'$ path } P \text{ in } G'' \text{ with robust value } < K+\tfrac{1}{2}.
\]
Hence the constructed Diff-RSP instance recovers the threshold structure of Most Secluded Path.

It remains only to verify the running time of the construction. For each node $v\in V$, the gadget $E_v$ contributes $4|V|$ chain edges and $4|V|+1$ chain vertices, together with one dummy vertex $d_v$, one entry vertex $c_v$, one dummy edge $\alpha_v$, and one entry edge $\beta_v$. Summing over all $v\in V$, the total number of gadget vertices and edges is therefore $O(|V|^2)$. The anchor edges contribute exactly $|E|$ additional edges, and the connector edges also contribute $|E|$ additional edges. Hence the graph $G''$ has size $O(|V|^2+|E|)$ and can be constructed from $G$ in $O(|V|^2+|E|)$ time. This completes the proof.


\subsection{Proofs for Section~\ref{sec:tsp}}


\subsubsection{Proof of Theorem~\ref{thm:diff-rtsp}}

The result follows immediately from Proposition~\ref{prop:tsp-short-infty} and Proposition~\ref{prop:tsp-global-budget}.


\subsubsection{Proof of Proposition~\ref{prop:tsp-short-infty}}

For each vertex $u\in V$, define
\[
    T_u:=\sum_{a\in E_{\mathrm{in}}(u)}\min\{\varepsilon,w_a\}.
\]
The weights $w^{\mathrm{wc}}$ can be computed in $O(|E|+|V|)$ time. Indeed, after initializing one value $T_u$ for each vertex, all quantities $T_u$ are computed by one pass over the edges. Then, for each edge $e=(i,u)$, we compute $\chi_e =\min\{\varepsilon,\,T_u-\min\{\varepsilon,w_e\}\}$ and set $w_e^{\mathrm{wc}}=w_e+\chi_e$ in $O(|E|)$ time. Fix now a Hamiltonian cycle $H\in\mathcal H$. For each vertex $u\in V$, let
\[
    e_H^-(u)\in H\cap E_{\mathrm{in}}(u),
    \qquad
    e_H^+(u)\in H\cap E_{\mathrm{out}}(u)
\]
be the unique tour edges entering and leaving $u$, respectively. We prove the claimed identity by matching upper and lower bounds.

First, take any feasible diffusion $(\Delta^+,\Delta^-)\in\mathcal D^{\mathrm{S},\infty}(\varepsilon)$. The perturbation of the tour cost can be written vertexwise as
\[
    \sum_{e\in H}(\Delta_e^+-\Delta_e^-)
    =
    \sum_{u\in V}
    \left(
        \Delta^+_{e_H^+(u)}
        -
        \Delta^-_{e_H^-(u)}
    \right).
\]
By flow conservation at vertex $u$,
\[
    \sum_{e\in E_{\mathrm{out}}(u)}\Delta_e^+
    =
    \sum_{a\in E_{\mathrm{in}}(u)}\Delta_a^-.
\]
Therefore,
\[
    \Delta^+_{e_H^+(u)}
    \le
    \sum_{e\in E_{\mathrm{out}}(u)}\Delta_e^+
    =
    \sum_{a\in E_{\mathrm{in}}(u)}\Delta_a^-,
\]
and hence
\[
    \Delta^+_{e_H^+(u)}
    -
    \Delta^-_{e_H^-(u)}
    \le
    \sum_{a\in E_{\mathrm{in}}(u)\setminus\{e_H^-(u)\}}\Delta_a^-.
\]
Under the short-term local budget, each incoming edge $a$ can lose at most
\[
    \Delta_a^- \le \min\{\varepsilon,w_a\}.
\]
Thus,
\[
    \Delta^+_{e_H^+(u)}
    -
    \Delta^-_{e_H^-(u)}
    \le
    \sum_{a\in E_{\mathrm{in}}(u)\setminus\{e_H^-(u)\}}
    \min\{\varepsilon,w_a\}.
\]
On the other hand, since $\Delta^+_{e_H^+(u)}\le \varepsilon$ and $\Delta^-_{e_H^-(u)}\ge 0$, we also have
\[
    \Delta^+_{e_H^+(u)}
    -
    \Delta^-_{e_H^-(u)}
    \le
    \varepsilon.
\]
Combining these two bounds gives
\[
    \Delta^+_{e_H^+(u)}
    -
    \Delta^-_{e_H^-(u)}
    \le
    \min\left\{
        \varepsilon,\,
        \sum_{a\in E_{\mathrm{in}}(u)\setminus\{e_H^-(u)\}}
        \min\{\varepsilon,w_a\}
    \right\}
    =
    \chi_{e_H^-(u)}.
\]
Summing over $u\in V$, we obtain
\[
    \sum_{e\in H}(\Delta_e^+-\Delta_e^-)
    \le
    \sum_{u\in V}\chi_{e_H^-(u)}
    =
    \sum_{e\in H}\chi_e.
\]
Since the feasible diffusion was arbitrary,
\[
    \max_{(\Delta^+,\Delta^-)\in\mathcal D^{\mathrm{S},\infty}(\varepsilon)}
    \sum_{e\in H}\bigl(w_e+\Delta_e^+-\Delta_e^-\bigr)
    \le
    \sum_{e\in H}(w_e+\chi_e)
    =
    \sum_{e\in H}w_e^{\mathrm{wc}}.
\]

It remains to prove that this upper bound is attainable. For each vertex $u\in V$, define
\[
    r_u
    :=
    \chi_{e_H^-(u)}
    =
    \min\left\{
        \varepsilon,\,
        T_u-\min\{\varepsilon,w_{e_H^-(u)}\}
    \right\}.
\]
By definition, $r_u$ is no larger than the total available capacity on the incoming non-tour edges into $u$. Hence we may choose numbers $\delta_a^u$, for $a\in E_{\mathrm{in}}(u)\setminus\{e_H^-(u)\}$, such that
\[
    0\le \delta_a^u\le \min\{\varepsilon,w_a\},
    \qquad
    \sum_{a\in E_{\mathrm{in}}(u)\setminus\{e_H^-(u)\}}\delta_a^u=r_u.
\]
For example, this can be done greedily over the incoming non-tour edges, since their total capacity is at least $r_u$. We now construct a diffusion $(\Delta_*^+,\Delta_*^-)$. For each vertex $u$, set
\[
    \Delta_{*,a}^-=\delta_a^u
    \quad
    \text{for all }
    a\in E_{\mathrm{in}}(u)\setminus\{e_H^-(u)\},
\]
and set
\[
    \Delta_{*,e_H^+(u)}^+=r_u.
\]
All remaining components of $\Delta_*^+$ and $\Delta_*^-$ are set to zero. 

We check feasibility. First, for every vertex $u$,
\[
    \sum_{e\in E_{\mathrm{out}}(u)}\Delta_{*,e}^+
    =
    \Delta_{*,e_H^+(u)}^+
    =
    r_u
    =
    \sum_{a\in E_{\mathrm{in}}(u)\setminus\{e_H^-(u)\}}\delta_a^u
    =
    \sum_{a\in E_{\mathrm{in}}(u)}\Delta_{*,a}^-.
\]
Thus the nodewise conservation constraints hold. Second, the short-term constraint holds because every edge with positive $\Delta_*^-$ satisfies
\[
    \Delta_{*,a}^-
    =
    \delta_a^u
    \le
    \min\{\varepsilon,w_a\}
    \le
    w_a.
\]
Finally, the local $\ell_\infty$ budget holds because
\[
    \Delta_{*,a}^-\le \varepsilon
    \quad\text{for every }a,
    \qquad
    \Delta_{*,e_H^+(u)}^+=r_u\le\varepsilon
    \quad\text{for every }u,
\]
and all other components are zero. Therefore
\[
    (\Delta_*^+,\Delta_*^-)\in\mathcal D^{\mathrm{S},\infty}(\varepsilon).
\]

For this feasible diffusion, no tour edge loses mass: indeed, $\Delta_{*,e_H^-(u)}^-=0$ for every $u\in V$. The only positive perturbation on a tour edge leaving $u$ is $\Delta_{*,e_H^+(u)}^+=r_u$. Therefore,
\[
    \sum_{e\in H}(\Delta_{*,e}^+-\Delta_{*,e}^-)
    =
    \sum_{u\in V}r_u
    =
    \sum_{u\in V}\chi_{e_H^-(u)}
    =
    \sum_{e\in H}\chi_e.
\]
Hence
\[
    \max_{(\Delta^+,\Delta^-)\in\mathcal D^{\mathrm{S},\infty}(\varepsilon)}
    \sum_{e\in H}\bigl(w_e+\Delta_e^+-\Delta_e^-\bigr)
    \ge
    \sum_{e\in H}(w_e+\chi_e)
    =
    \sum_{e\in H}w_e^{\mathrm{wc}}.
\]
Together with the upper bound, this proves the fixed-tour identity.

Since the identity holds for every Hamiltonian cycle $H$, minimizing the robust objective under $\mathcal D^{\mathrm{S},\infty}(\varepsilon)$ is equivalent to solving the ordinary TSP instance with edge weights $w^{\mathrm{wc}}$. The weights $w^{\mathrm{wc}}$ are computable in $O(|E|+|V|)$ time, so the construction of the ordinary TSP instance has the claimed complexity.


\subsubsection{Proof of Proposition~\ref{prop:tsp-global-budget}}

Fix a Hamiltonian cycle $H\in\mathcal H$, and set
\[
    r:=\min\left\{\frac{\varepsilon}{2},\,S-W(H)\right\}.
\]
We first prove the upper bound. Let $(\Delta^+,\Delta^-)\in \mathcal D(\varepsilon)$, where $\mathcal D(\varepsilon)$ is either $\mathcal D^{\mathrm{S},1}(\varepsilon)$ or $\mathcal D^{\mathrm{L},1}(\varepsilon)$. By conservation,
\[
    \sum_{e\in E}\Delta_e^+
    =
    \sum_{e\in E}\Delta_e^-.
\]
Since $\|(\Delta^+,\Delta^-)\|_1\le \varepsilon$, it follows that
\[
    \sum_{e\in E}\Delta_e^+
    =
    \sum_{e\in E}\Delta_e^-
    \le
    \frac{\varepsilon}{2}.
\]
Therefore,
\[
    \sum_{e\in H}(\Delta_e^+-\Delta_e^-)
    \le
    \sum_{e\in H}\Delta_e^+
    \le
    \sum_{e\in E}\Delta_e^+
    \le
    \frac{\varepsilon}{2}.
\]
Thus the cost of $H$ is at most $W(H)+\varepsilon/2$.

We also have a second upper bound. In both the short-term and long-term regimes, the post-diffusion weights are nonnegative: in the short-term case this follows from $\Delta^-\le w$, and in the long-term case from $\Delta^-\le w+\Delta^+$. Hence
\[
    \sum_{e\in H}\bigl(w_e+\Delta_e^+-\Delta_e^-\bigr)
    \le
    \sum_{e\in E}\bigl(w_e+\Delta_e^+-\Delta_e^-\bigr).
\]
By conservation, the total post-diffusion mass is
\[
    \sum_{e\in E}\bigl(w_e+\Delta_e^+-\Delta_e^-\bigr)
    =
    \sum_{e\in E}w_e
    =
    S.
\]
Combining the two upper bounds gives
\[
    \sum_{e\in H}\bigl(w_e+\Delta_e^+-\Delta_e^-\bigr)
    \le
    \min\left\{
        W(H)+\frac{\varepsilon}{2},
        S
    \right\}.
\]
Since $(\Delta^+,\Delta^-)$ was arbitrary, this proves the desired upper bound on the inner maximization.

It remains to show that the bound is attainable. Since
\[
    r\le S-W(H)=\sum_{e\notin H}w_e,
\]
we may choose numbers $\delta_e$, for $e\notin H$, such that
\[
    0\le \delta_e\le w_e,
    \qquad
    \sum_{e\notin H}\delta_e=r.
\]
For example, this can be done greedily over the edges outside $H$. For each vertex $u\in V$, let
\[
    e_H^+(u)\in H\cap E_{\mathrm{out}}(u)
\]
denote the unique tour edge leaving $u$, and define
\[
    q_u:=\sum_{a\in E_{\mathrm{in}}(u)\setminus H}\delta_a.
\]
We construct a diffusion $(\Delta_*^+,\Delta_*^-)$. For each off-tour edge $e\notin H$, set
\[
    \Delta_{*,e}^-:=\delta_e,
\]
and for each vertex $u\in V$, set
\[
    \Delta_{*,e_H^+(u)}^+:=q_u.
\]
All remaining components of $\Delta_*^+$ and $\Delta_*^-$ are set to zero.

We check feasibility. For every vertex $u\in V$,
\[
    \sum_{e\in E_{\mathrm{out}}(u)}\Delta_{*,e}^+
    =
    \Delta_{*,e_H^+(u)}^+
    =
    q_u
    =
    \sum_{a\in E_{\mathrm{in}}(u)\setminus H}\delta_a
    =
    \sum_{a\in E_{\mathrm{in}}(u)}\Delta_{*,a}^-.
\]
Thus the nodewise conservation constraints hold. Moreover, $\Delta_*^-\le w$ because each positive component satisfies $\Delta_{*,e}^-=\delta_e\le w_e$. Hence the short-term constraint holds. Since $\Delta_*^+\ge 0$, we also have $\Delta_*^-\le w+\Delta_*^+$, so the long-term constraint holds as well. Finally,
\[
    \|(\Delta_*^+,\Delta_*^-)\|_1
    =
    \sum_{e\in E}\Delta_{*,e}^+
    +
    \sum_{e\in E}\Delta_{*,e}^-
    =
    \sum_{u\in V}q_u+\sum_{e\notin H}\delta_e
    =
    2r
    \le
    \varepsilon.
\]
Thus $(\Delta_*^+,\Delta_*^-)$ is feasible for both
$\mathcal D^{\mathrm{S},1}(\varepsilon)$ and
$\mathcal D^{\mathrm{L},1}(\varepsilon)$.

For this feasible diffusion, no tour edge loses mass, because $\Delta_{*,e}^-=0$ for every $e\in H$. The only perturbation on tour edges is the mass injected into the unique outgoing tour edge at each vertex. Therefore,
\[
    \sum_{e\in H}(\Delta_{*,e}^+-\Delta_{*,e}^-)
    =
    \sum_{u\in V}q_u
    =
    \sum_{e\notin H}\delta_e
    =
    r.
\]
Consequently,
\[
    \sum_{e\in H}\bigl(w_e+\Delta_{*,e}^+-\Delta_{*,e}^-\bigr)
    =
    W(H)+r
    =
    \min\left\{
        W(H)+\frac{\varepsilon}{2},
        S
    \right\}.
\]
This proves the fixed-tour identity.

It remains only to derive the reduction statement. Since the function
\[
    x\mapsto \min\left\{x+\frac{\varepsilon}{2},\,S\right\}
\]
is nondecreasing in $x$, minimizing the robust value over $H\in\mathcal H$ is equivalent to minimizing $W(H)$. Therefore
\[
    \mathrm{OPT}_{\mathrm{RTSP}}(w,\mathcal D(\varepsilon))
    =
    \min\left\{
        \mathrm{OPT}(w)+\frac{\varepsilon}{2},
        S
    \right\}.
\]
Thus, under either $\mathcal D^{\mathrm{S},1}(\varepsilon)$ or $\mathcal D^{\mathrm{L},1}(\varepsilon)$, Diff-RTSP reduces to ordinary TSP with the original weights $w$, followed by scalar postprocessing.


\subsubsection{Proof of Lemma~\ref{lem:tsp-long-infty-upper}}

Fix a Hamiltonian cycle $H\in\mathcal H$ and let $(\Delta^+,\Delta^-)\in\mathcal D^{\mathrm{L},\infty}(\varepsilon)$ be arbitrary. We prove two upper bounds on the perturbation term $\sum_{e\in H}(\Delta_e^+-\Delta_e^-)$.

First, since $H$ contains exactly $n$ edges and the local budget gives $\Delta_e^+\le \varepsilon$ for every $e\in E$, we have
\[
    \sum_{e\in H}(\Delta_e^+-\Delta_e^-)
    \le
    \sum_{e\in H}\Delta_e^+
    \le
    n\varepsilon.
\]
Therefore,
\[
    \sum_{e\in H}\bigl(w_e+\Delta_e^+-\Delta_e^-\bigr)
    \le
    W(H)+n\varepsilon.
\]

Second, conservation implies that the total perturbation mass injected into the graph equals the total perturbation mass drained from the graph:
\[
    \sum_{e\in E}\Delta_e^+
    =
    \sum_{e\in E}\Delta_e^-.
\]
Hence
\[
    \sum_{e\in H}(\Delta_e^+-\Delta_e^-)
    =
    -\sum_{e\notin H}(\Delta_e^+-\Delta_e^-)
    =
    \sum_{e\notin H}(\Delta_e^- - \Delta_e^+).
\]
For each edge $e\notin H$, the local budget gives
\[
    \Delta_e^- - \Delta_e^+
    \le
    \Delta_e^-
    \le
    \varepsilon.
\]
Moreover, the long-term constraint $\Delta^- \le w+\Delta^+$ gives
\[
    \Delta_e^- - \Delta_e^+
    \le
    w_e.
\]
Combining these two inequalities yields
\[
    \Delta_e^- - \Delta_e^+
    \le
    \min\{\varepsilon,w_e\}
    =
    c_e.
\]
Summing over $e\notin H$, we obtain
\[
    \sum_{e\in H}(\Delta_e^+-\Delta_e^-)
    =\sum_{e\notin H}(\Delta^-_e-\Delta_e^+)
    \le
    \sum_{e\notin H}c_e
    =
    C-\sum_{e\in H}c_e.
\]
Thus,
\[
    \sum_{e\in H}\bigl(w_e+\Delta_e^+-\Delta_e^-\bigr)
    \le
    W(H)+C-\sum_{e\in H}c_e
    =
    C+\sum_{e\in H}(w_e-c_e).
\]
Since the feasible diffusion $(\Delta^+,\Delta^-)$ was arbitrary, the two upper bounds imply
\[
    \max_{(\Delta^+,\Delta^-)\in\mathcal D^{\mathrm{L},\infty}(\varepsilon)}
    \sum_{e\in H}\bigl(w_e+\Delta_e^+-\Delta_e^-\bigr)
    \le
    \min\left\{
        W(H)+n\varepsilon,\,
        C+\sum_{e\in H}(w_e-c_e)
    \right\}.
\]
It remains to derive the bound on the optimal robust value. Taking the minimum over $H\in\mathcal H$ in the fixed-tour upper bound gives
\[
\begin{aligned}
    \mathrm{OPT}_{\mathrm{RTSP}}\bigl(w,\mathcal D^{\mathrm{L},\infty}(\varepsilon)\bigr)
    &\le
    \min_{H\in\mathcal H}
    \min\left\{
        W(H)+n\varepsilon,\,
        C+\sum_{e\in H}(w_e-c_e)
    \right\} \\
    &=
    \min\left\{
        \min_{H\in\mathcal H}\bigl(W(H)+n\varepsilon\bigr),\,
        \min_{H\in\mathcal H}\left(C+\sum_{e\in H}(w_e-c_e)\right)
    \right\} \\
    &=
    \min\left\{
        \mathrm{OPT}(w)+n\varepsilon,\,
        C+\mathrm{OPT}(w-c)
    \right\}.
\end{aligned}
\]
Finally, the quantities $c$, $w-c$, and $C$ are computed by one pass over the edges, which takes $O(|E|)$ time.

\end{document}